\documentclass[11pt]{article}

\usepackage{amsmath}  % Nodig voor het gebruik van \numberwithin
\usepackage{amsthm} % nodig voor het gebruik van \begin{proof}
\usepackage{amssymb} % nodig voor \Bbb en \mathbf (wat beide hetzelfde is, maar \Bbb wordt afgeraden)

\usepackage{enumitem}
\usepackage[backref=page]{hyperref}  %

\usepackage{color} % Voor het gebruik van kleuren

\newtheorem{thm}{Theorem}
\newtheorem{inspr}[thm]{}

\newenvironment{definitie}{\begin{itemize}\item[ ]\hspace{-26pt}\bf Definition \rm }{\end{itemize}}
\newenvironment{notatie}{\begin{itemize}\item[ ]\hspace{-26pt}\bf Notation \rm }{\end{itemize}}
\newenvironment{voorbeeld}{\begin{itemize}\item[ ]\hspace{-26pt}\bf Example \rm }{\end{itemize}}
\newenvironment{stelling}{\begin{itemize}\item[ ]\hspace{-26pt}\bf Theorem \rm }{\end{itemize}}
\newenvironment{propositie}{\begin{itemize}\item[ ]\hspace{-26pt}\bf Proposition \rm }{\end{itemize}}
\newenvironment{lemma}{\begin{itemize}\item[ ]\hspace{-26pt}\bf Lemma \rm }{\end{itemize}}
\newenvironment{opmerking}{\begin{itemize}\item[ ]\hspace{-26pt}\bf Remark \rm }{\end{itemize}}
\newenvironment{voorwaarde}{\begin{itemize}\item[ ]\hspace{-26pt}\bf Condition \rm }{\end{itemize}}
\newenvironment{probleem}{\begin{itemize}\item[ ]\hspace{-26pt}\bf Problem \rm }{\end{itemize}}
\newenvironment{gevolg}{\begin{itemize}\item[ ]\hspace{-26pt}\bf Corollary \rm }{\end{itemize}}
\newenvironment{niets}{\begin{itemize}\item[ ]\hspace{-26pt}\bf   \rm }{\end{itemize}}

\newcommand{\defin}{\begin{inspr}\begin{definitie}}  %\def already defined
\newcommand{\edefin}{\end{definitie}\end{inspr}}

\newcommand{\notat}{\begin{inspr}\begin{notatie}}  %\not already defined
\newcommand{\enotat}{\end{notatie}\end{inspr}}

\newcommand{\voorb}{\begin{inspr}\begin{voorbeeld}}  %\not already defined
\newcommand{\evoorb}{\end{voorbeeld}\end{inspr}}

\newcommand{\stel}{\begin{inspr}\begin{stelling}}
\newcommand{\estel}{\end{stelling}\end{inspr}}

\newcommand{\prop}{\begin{inspr}\begin{propositie}}
\newcommand{\eprop}{\end{propositie}\end{inspr}}

\newcommand{\lem}{\begin{inspr}\begin{lemma}}
\newcommand{\elem}{\end{lemma}\end{inspr}}

\newcommand{\opm}{\begin{inspr}\begin{opmerking}}
\newcommand{\eopm}{\end{opmerking}\end{inspr}}

\newcommand{\voorw}{\begin{inspr}\begin{voorwaarde}}
\newcommand{\evoorw}{\end{voorwaarde}\end{inspr}}

\newcommand{\probl}{\begin{inspr}\begin{probleem}}
\newcommand{\eprobl}{\end{probleem}\end{inspr}}

\newcommand{\gev}{\begin{inspr}\begin{gevolg}}
\newcommand{\egev}{\end{gevolg}\end{inspr}}

\newcommand{\nul}{\begin{inspr}\begin{niets}}
\newcommand{\enul}{\end{niets}\end{inspr}}

\newcommand{\bew}{\vspace{-0.3cm}\begin{itemize}\item[ ] \bf Proof\rm: }
\newcommand{\ebew}{\hfill $\qed$ \end{itemize}}

\newcommand{\ssnl}{\vskip 3pt} % Het is noodzakelijk dat er voor de instructie vspace 
\newcommand{\snl}{\vskip 7pt} % Het is noodzakelijk dat er voor de instructie vspace een lege lijn staat
\newcommand{\nl}{\vskip 12pt} % Kunnen we dit hierin opnemen?

\newcommand{\ot}{\otimes}

\newcommand{\inv}{^{-1}}

\newcommand{\tussenen}{\qquad\quad\text{and}\qquad\quad}

\newcommand{\rood}{\color{red}}

\numberwithin{thm}{section}   % Zorgt ervoor dat de nummering bestaat uit ?.?
\numberwithin{equation}{section} % Zorgt ervoor dat de nummering bestaat uit ?.?

\newcommand{\keepcomment}[1]{}
\newcommand{\oldcomment}[1]{}

\begin{document}

%%%%%%%%%%%%%%%%%%%%%%%%%%%%%%%%%%%%%%%%%%%%%%%%%%%%%%%%%%%%%%%%%%%%%%
%
% Bestand artikela.tex
%
%%%%%%%%%%%%%%%%%%%%%%%%%%%%%%%%%%%%%%%%%%%%%%%%%%%%%%%%%%%%%%%%%%%%%

\centerline{\bf \Large Reflections on the Larson-Sweedler theorem}
\centerline{\bf\large for (weak) multiplier Hopf algebras} 
\vspace{13 pt}
\centerline{\it A.\ Van Daele \rm ($^*$)}
\vspace{20 pt}
{\bf Abstract} 
\nl 
Let $A$ be an algebra with identity and  $\Delta:A\to A\ot A$  a coproduct that admits a counit. If there exist a faithful left integral and a faithful right integral, one can construct an antipode and $(A,\Delta)$ is a Hopf algebra. This is the Larson-Sweedler theorem.
\ssnl
There are generalizations of this result for multiplier Hopf algebras, weak Hopf algebras and weak multiplier Hopf algebras. In the case of a multiplier Hopf algebra, the existence of a counit can be weakened and can be replaced by the requirement that the coproduct is full. A similar result is true for weak multiplier Hopf algebras. 
\ssnl
What we show in this note is that in fact the result for multiplier Hopf algebras can still be obtained without the condition of fullness of the coproduct. As it turns out, this property will already follow from the other conditions. Consequently, also in the original theorem for Hopf algebras, the existence of a counit is a consequence of the other conditions. This slightly generalizes the original result.
\ssnl
The situation for weak multiplier Hopf algebras seems to be more subtle. We discuss the problems and see what is still possible here.
\ssnl
We consider these results in connection with the development of the theory of locally compact quantum groups. This is discussed in an appendix.
\nl \nl
Date: {\it 6 May 2024} 

\vskip 6cm

\hrule
\vskip 7 pt
\begin{itemize}
\item[($^*$)] Department of Mathematics, KU Leuven, Celestijnenlaan 200B,\newline
B-3001 Heverlee (Belgium). E-mail: {\tt alfons.vandaele@kuleuven.be}
\end{itemize}

\newpage

%%%%%%%%%%%%%%%%%%%%%%%%%%%%%%%%%%%%%%%%%%%%%%%%%%%%%
%
% Bestand: artikel0.tex
%
%%%%%%%%%%%%%%%%%%%%%%%%%%%%%%%%%%%%%%%%%%%%%%%%%%%%%%

\setcounter{section}{-1}  % Dit zorgt ervoor dat we met 0 beginnen voor de inleiding

\section{\hspace{-17pt}. Introduction} \label{s:intro} % \input artikel0.tex

In this paper we work with an algebra $A$ over the complex numbers. We do not assume that it has an identity, but we require the product to be non-degenerate (as a bilinear map). Remark that this is automatic when the algebra is unital.
\ssnl
 We consider the algebra $A\ot A$. It is again a non-degenerate algebra. We use $M(A)$ and $M(A\ot A)$ for the multiplier algebras of $A$ and $A\ot A$ respectively. 
\ssnl
We consider the obvious inclusions
$$A\ot A\subseteq M(A)\ot M(A) \subseteq M(A\ot A).$$ 
Observe that these inclusions, in general, are strict.
We refer to the item {\it Notations, conventions and basic references} for details.
\snl
In general, a coproduct on a vector space $V$  is defined as a linear map $\Delta:V\to V\ot V$ satisfying coassociativity $(\Delta\ot\iota)\Delta=(\iota\ot\Delta)\Delta$. We use $\iota$ for the identity map from $V$ to itself. This is the usual definition when $V$ is a unital algebra. In the more general case of a non-degenerate algebra $ A$, possibly without identity, it would be far too restrictive for a coproduct to assume that it takes values in $A\ot A$. But then we have a problem with the formulation of coassociativity in its usual form as the maps $\Delta\ot\iota$ and $\iota\ot\Delta$ are not defined on the range of $\Delta$. 
\ssnl
A correct description of coassociativity is crucial for the treatment of the Larson-Sweedler theorem in this note. 
\ssnl
This problem is discussed in detail in a note we have written on the coproduct on non-unital algebras, see \cite{VD-refl}. We will briefly recall some aspects of this in a preliminary section of this paper (see Section \ref{s:prel}).
\nl
\bf Content of the paper\rm
\nl
In {\it Section} \ref{s:prel}, we first recall the different regularity conditions of a coproduct, needed to formulate coassociativity. Further we mention some other properties of a coproduct that are important for the rest of the paper. There is the notion of {\it fullness} of a coproduct and the possible {\it existence of a counit}. The two properties are related and this relation is discussed. For the main results in this note, we will need regularity of some of the canonical maps, while fullness of the coproduct and the existence of a counit will follow from the assumptions. \oldcomment{\rood Verify and compare with the other paper.}
\ssnl
We will need to recall other notions, like an antipode, left and right integrals, ... . This will be done further in the paper, at those places where these concepts first appear.
\ssnl
In {\it Section} \ref{s:ls-mha} we discuss the {\it main result} for {\it multiplier Hopf algebras}. We start with a pair $(A,\Delta)$ of a non-degenerate algebra $A$ and a coproduct $\Delta$ on $A$. We assume the existence of a  left integral and a right integral. Under the appropriate  regularity conditions of $\Delta$ and faithfulness properties of the integrals,  we show that there is an antipode and that actually $(A,\Delta)$ is a  multiplier Hopf algebra. This is the Larson-Sweedler theorem for multiplier Hopf algebras. In an earlier paper on the subject (\cite{VD-W-ls}), we have proven the same result but under the  assumption that the coproduct is regular and that it is full. Here we show that fullness is not needed to prove the result but that it actually follows from the other assumptions. 
\ssnl
In {\it Section} \ref{s:ls-wmha} we consider the more general case of a {\it weak} multiplier Hopf algebra. The starting point is the same, but now we use a more general notion of integrals. We have to assume the existence of a faithful set of left integrals and a faithful set of right integrals. Again the existence of an antipode is proven under these conditions and the pair $(A,\Delta)$ now turns out to be a regular weak  multiplier Hopf algebra. This is the Larson-Sweedler theorem for weak multiplier Hopf algebras. 
\ssnl
In an earlier work on the subject (\cite{K-VD}), fullness of the coproduct is used to obtain the result. In this note, we attempt to avoid this assumption and work along the lines of the arguments we developed in the previous section in the case of a multiplier Hopf algebra. 
\ssnl
In {\it Section} \ref{s:concl} we draw some conclusions and consider some possible further research.
\ssnl
We have included two appendices. In {\it Appendix} \ref{sect:appA} we treat an aspect of integrals on multiplier Hopf algebras that is not important for the theory of multiplier Hopf algebras as such, but rather for the Larson-Sweedler theorem as we treat it in this paper. In {\it Appendix} \ref{sect:appB} we discuss the relevance of the Larson-Sweedler theorem as we prove it here for the understanding of the basic ideas behind the  definitions of \emph{locally compact quantum groups}, both in the C$^*$-algebraic setting as well as in the von Neumann algebraic framework. In fact, it is the work on locally compact quantum groups that inspired us to consider this more advanced approach to the Larson-Sweedler theorem.
\nl
\bf Notations, conventions and basic references\rm
\nl
As mentioned already, we only work with algebras over the complex numbers, although we believe that most of the results are true for algebras over other fields. 
\ssnl
If the algebra has an identity, we use $1$ to denote it. If it has no identity, we do assume that the product is non-degenerate. For such an algebra $A$ we denote by $M(A)$ the multiplier algebra. The identity in $M(A)$ is also denoted by $1$. 
\ssnl
The tensor product $A\ot A$ is again non-degenerate and we use $M(A\ot A)$ for its multiplier algebra. 
\ssnl
The space of linear functionals on $A$ is denoted by $A'$. A linear functional $\omega$ on $A$ is called {\it faithful} if the linear maps $a\mapsto \omega(a\,\cdot\,)$ and $a\mapsto \omega(\,\cdot\,a)$ are both injective from $A$ to $A'$. We will also consider {\it single sided} faithful linear functionals (see Definition \ref{defin:fff} in Section \ref{s:ls-mha}), as well as faithful sets of linear functionals (see Definition \ref{defin:3.11c}). 
\ssnl
If $A$ is finite-dimensional, a linear functional  is automatically faithful if it is either left or right faithful. This is no longer true in the infinite-dimensional case.
\ssnl
A coproduct on the algebra $A$ will always be a linear map from $A$ to the multiplier algebra $M(A\ot A)$. It is assumed to satisfy the conditions formulated in Section \ref{s:prel}. 
\ssnl
We use the {\it leg numbering notation} for a coproduct $\Delta$. We write $\Delta_{12}(a)$, where $a\in A$, for $\Delta(a)\ot 1$ as sitting in
$M(A\ot A\ot A)$. Similarly, we use $\Delta_{23}(a)$ for $1\ot \Delta(a)$ in $M(A\ot A\ot A)$. Finally, $\Delta_{13}(a)$ is as $\Delta(a)$ sitting with its first leg in the first factor and with its second leg in the third factor of $M(A\ot A\ot A)$. More precisely it  is $(\zeta\ot\iota)\Delta_{23}(a)$ where $\zeta$ is the flip on $A\ot A$ and $\zeta\ot\iota$ is extended to $M(A\ot A\ot A)$.
\ssnl
We will sometimes use the Sweedler notation for the coproduct. This is justified provided  we have the correct coverings. This is explained in \cite{D-VD}. See also \cite{VD-tools} and the more recent work on this topic, found in \cite{VD-sw}. Remark however that the Sweedler notation is just what it says, a notation, and nothing more. It is simply another way of writing down results and formulas to obtain the results. It has the advantage of making formulas more transparent but the disadvantage that it is less rigorous. Still, it should always be clear how to make the arguments completely rigorous by avoiding the use of the Sweedler notation. We illustrate this point of view by providing two proofs of similar results, one using the Sweedler notation and another one where the Sweedler notation is not used, see Proposition \ref{prop:rTl}. 
\ssnl
If an antipode exists, we will denote it by $S$. The defining properties for an antipode are more complicated, certainly in the case of a weak multiplier Hopf algebra. We will discuss these when we meet the antipode in the course of the paper.
\ssnl
For the theory of Hopf algebras we have the standard references \cite{A} and \cite{S}, see also the more recent work \cite{R}.
For the theory of multiplier Hopf algebras, the main (original) reference is \cite{VD-mha} and for the theory of multiplier Hopf algebras with integrals, sometimes called algebraic quantum groups, the main reference is \cite{VD-afw}. Weak multiplier Hopf algebras are studied in a number of papers. See \cite{VD-W1, VD-W2} and \cite{VD-W3} (and also \cite{VD-W0}). 
\nl
\bf Acknowledgments \rm
\nl
I am grateful for having the opportunity to continue my research at the KU Leuven after my retirement. Further I  like to thank my colleagues and friends at the Department of Mathematics of the University of Oslo and Trondheim (Norway), where part of this work has been written, for their hospitality and the nice working atmosphere.

%%%%%%%%%%%%%%%%%%%%%%%%%%%%%%%%%%%%%%%%%%%%%%%%%%%%%
%
% Section: Preliminaries
%
%%%%%%%%%%%%%%%%%%%%%%%%%%%%%%%%%%%%%%%%%%%%%%%%%%%%%%

\section{\hspace{-17pt}. Preliminaries} \label{s:prel} % \input artikel1.tex

As mentioned already in the introduction, in this paper we work with an algebra $A$ over the complex numbers, it is not assumed to be unital, but we do require that the product is non-degenerate. The tensor product algebra $A\ot A$ is again a non-degenerate algebra. We use $M(A)$ and $M(A\ot A)$ for the multiplier algebras of $A$ and $A\ot A$ respectively. 
\ssnl
For a coproduct $\Delta$  on a non-unital algebra, it is too restrictive to assume that it takes values in $A\ot A$. We need to consider linear maps from $A$ into the multiplier algebra $M(A\ot A)$ of $A\ot A$. This needs some extra considerations that are not relevant if the algebra has an identity.
\nl
\bf Coproducts on non-degenerate algebras \rm
\nl
For an algebra $A$, not necessarily unital, a coproduct is  a linear map $\Delta:A\to M(A\ot A)$ to begin with. We cannot simply express coassociativity in its usual form $(\Delta\ot\iota)\Delta=(\iota\ot\Delta)\Delta$ as the maps $\Delta\ot\iota$ and $\iota\ot\Delta$ are not defined on the range of $\Delta$. 
The problem can be overcome in different ways. 
\ssnl
One possibility is to assume that $\Delta$ is a \emph{homomorphism} from $A$ to $M(A\ot A)$ and that it is \emph{non-degenerate}. That means that the subspaces of $A\ot A$, spanned by elements of the form
\begin{equation*}
\Delta(a)(c\ot b)
\tussenen
(c\ot b)\Delta(a)
\end{equation*}
are both all of $A\ot A$. Under this condition, it is possible to extend the homomorphisms $\Delta\ot\iota$ and $\iota\ot\Delta$ in a unique way to unital homomorphisms from $M(A\ot A)$ to  $M(A\ot A\ot A)$. We usually denote these extensions by the same symbols and then coassociativity is still formulated as $(\Delta\ot\iota)\Delta=(\iota\ot\Delta)\Delta$.
\ssnl
Unfortunately, non-degeneracy of the coproduct is too restrictive in the case of weak multiplier Hopf algebras. The condition can be weakened so that still these homomorphisms can be extended and so that coassociativity can  be formulated in its usual form.
\ssnl
All of this is explained in detail, in the recent note on coproducts for non-unital algebras, see \cite{VD-refl}.
\ssnl
This  is one way to solve the problem. 

In this note however, we will use a second possibility (which is in fact the more common one). It is also discussed in detail in \cite{VD-refl}, where moreover the relation between the two approaches is given.
\ssnl
We recall the basic concepts here. The starting point is just a linear map $A\mapsto M(A\ot A)$. We associate the following four maps.

\notat\label{notat:Tmaps}
Given a linear map $\Delta:A\to M(A\ot A)$, we consider the  maps from $A\ot A$ to $M(A\ot A)$, defined by
\begin{align*}
T_1(a\ot b)&=\Delta(a)(1\ot b)
\tussenen
T_2 (c\ot a)=(c\ot 1)\Delta(a) \\
T_3(a\ot b)&=(1\ot b)\Delta(a)
\tussenen
T_4(c\ot a)=\Delta(a)(c\ot 1)
\end{align*}
where $a,b,c\in A$.
\enotat

We call these maps the {\it canonical maps} associated to $\Delta$. If a canonical map has range in $A\ot A$, we say that it is {\it regular}. If all four of the canonical maps have range in $A\ot A$, we call $\Delta$ regular. 
\ssnl
The assumptions on the ranges of these canonical maps are meaningless for an algebra with identity as then these conditions are trivially satisfied. The same is true for the notion of regularity of a coproduct.
\ssnl
We formulate coassociativity of such a linear map $\Delta$ in various situations.

\defin\label{defin:coass} 
i) If the maps $T_1$ and $T_2$ are regular, we say that $\Delta$ is coassociative if 
\begin{equation}
(c\ot 1\ot 1)(\Delta\ot\iota)(\Delta(a)(1\ot b))
=(\iota\ot\Delta)((c\ot 1)\Delta(a))(1\ot 1\ot b)\label{eqn:1.2a}
\end{equation}
for all $a,b,c\in A$. 
\ssnl
ii) If the maps $T_3$ and $T_4$ are regular, we say that $\Delta$ is coassociative if
\begin{equation}
(\Delta\ot\iota)((1\ot b)\Delta(a))(c\ot 1\ot 1)
=(1\ot 1\ot b)(\iota\ot\Delta)(\Delta(a)(c\ot 1))\label{eqn:1.2b}
\end{equation}
for all $a,b,c\in A$. 
\edefin

Remark that regularity of the maps $T_1$ and $T_2$ is needed  for Equation (\ref{eqn:1.2a}) to make sense. Similarly we need the regularity of $T_3$ and $T_4$ for Equation \ref{eqn:1.2b}. One form is equivalent with the other when we replace $A$ by $A^\text{op}$ or $\Delta$ by $\Delta^\text{cop}$.
\ssnl
Equation (\ref{eqn:1.2a}) can be reformulated as 
$(T_2\ot\iota)(\iota\ot T_1)=(\iota\ot T_1)(T_2\ot\iota)$. Similar alternative expressions are possible for the other case. 
\ssnl
The first condition is the more common one.
However, apart from these two, there still are  other possibilities.  

\defin\label{defin:1.11a} 
i) Assume that $T_1$ and $T_4$ are regular. Then we call $\Delta$ coassociative if
\begin{equation}
((\Delta\ot\iota)(\Delta(a)(1\ot b)))(c\ot 1\ot 1)
=((\iota\ot\Delta)(\Delta(a)(c\ot 1)))(1\ot 1\ot b)\label{eqn:1.3a}
\end{equation}
for all $a,b,c\in A$.
\ssnl
ii) \oldcomment{This is wrong!}{}
Assume that $T_2$ and $T_3$ are regular. Then we call $\Delta$ coassociative if
\begin{equation}
(c\ot 1\ot 1)((\Delta\ot\iota)((1\ot b)\Delta(a)))
=(1\ot 1\ot b)((\iota\ot\Delta)((c\ot 1)\Delta(a)))\label{eqn:1.3b}
\end{equation}
for all $a,b,c\in A$.
\edefin

Remark that if the four canonical maps are regular, all these conditions are equivalent. This is shown by using the non-degeneracy of the product in $A$.
\snl
This takes us to the following basic definition, used in this paper.

\defin\label{defin:copr}
Let $A$ be a non-degenerate algebra. A coproduct on $A$ is a \emph{homomorphism} $\Delta$ from  $A$ to $M(A\ot A)$ satisfying \emph{coassociativity} (in any of these forms).
\edefin

So we only call a linear map $\Delta:A\to M(A\ot A)$ a coproduct if either $T_1$ and $T_2$, $T_3$ and $T_4$, $T_1$ and $T_4$ or $T_2$ and $T_3$ have range in $A\ot A$. These are precisely the four cases that allow the formulation of coassociativity.
\ssnl
Remark that in the literature on multiplier Hopf algebras and weak multiplier Hopf algebras, for non-regular coproducts,  only the first case, where $T_1$ and $T_2$ are assumed to map in $A\ot A$, have been considered. For regular multiplier Hopf algebroids, as introduced in \cite{T-VD}, the cases with $T_1$ combined with $T_4$ and $T_2$ combined with $T_3$ are used. These combinations are also used in the more recent paper on single-sided multiplier Hopf algebras \cite{VD-ssmhas}.
\ssnl
The definition, as formulated here, is symmetric in the following sense. If $\Delta$ is a coproduct on $A$, then it is also a coproduct on the opposite algebra $A^{\text{op}}$.  Also  the flipped map $\Delta{^\text{cop}}$ is a coproduct on $A$. These two cases give a regular coproduct if there is regularity in the original situation.
\ssnl
There are also ways to express coassociativity if only one canonical map is regular (see \cite{VD-refl} and also \cite{VD-ssmhas} , but this situation does not occur in this paper.
\ssnl
Finally remark that for several notions and results that follow, it is not always necessary to have that the linear map $\Delta:A\to M(A\ot A)$ is a homomorphism, neither that it is coassociative. We will mention this when appropriate. 
\nl
\bf The legs of a coproduct\rm
\nl
We recall the following definition, see e.g.\ Section 1 in \cite{VD-W-ls}.

\defin\label{defin:full-cp}
Let $\Delta$ be a linear map from $A$ to $M(A\ot A)$. Assume that the canonical maps $T_1$ and $T_2$ are regular. Denote by $V$ and $W$ the smallest subspaces of $A$ satisfying
\begin{equation}
\Delta(a)(1\ot b)\in V\ot A
\qquad\quad\text{and}\qquad\quad
(c\ot 1)\Delta(a)\in A\ot W
\end{equation}
for all $a,b,c\in A$. We will call $V$ the {\it left leg} of $\Delta$ and $W$ the {\it right leg} of $\Delta$. A coproduct is called {\it full} if $V=A$ and $W=A$.
\edefin

We have the obvious definitions in the other cases. Moreover, it is an easy consequence of the non-degeneracy of the product in $A$ that e.g.\ if both the canonical maps $T_1$ and $T_3$ are regular, then the smallest subspaces $V$ and $V'$ satisfying 
\begin{equation}
\Delta(a)(1\ot b)\in V\ot A
\qquad\quad\text{and}\qquad\quad
(1\ot b)\Delta(a)\in V'\ot A
\end{equation}
for all $a,b\in A$, will be the same. Similarly for the right leg if $T_2$ and $T_4$ are regular. \oldcomment{Is this still relevant?}
\ssnl
We can characterize these legs of $\Delta$ as follows, see Propositions 1.5 and 1.6 in \cite{VD-W-ls}. \oldcomment{\rood Perhaps we do not need to include this proof here, certainly not in the short version.}{}

\prop \label{prop:1.6c}
i) If $T_1$ is regular then the left leg $V$ of $\Delta$ is equal to the linear span of elements of the form $(\iota\ot \omega)(\Delta(a)(1\ot b))$, with $a,b\in A$ and $\omega\in A'$.  If $T_3$ is regular, it is the span of the elements $(\iota\ot\omega)((1\ot b)\Delta(a))$.
\ssnl
ii) If $T_2$ is regular then the right leg $W$ of $\Delta$ is equal to the linear span of elements of the form $(\omega\ot\iota)((c\ot 1)\Delta(a))$, with $a,c\in A$ and $\omega\in A'$.  If $T_4$ is regular, it is the span of the elements $(\omega\ot\iota)(\Delta(a)(c\ot 1))$.
\eprop
\bew
Assume that $T_1$ is regular. Assume that $V$ is a subspace of $A$ such that $\Delta(a)(1\ot b)\subseteq V\ot A$ for all $a,b$. Denote by $V'$ the space spanned by elements of the form $(\iota\ot \omega)(\Delta(a)(1\ot b))$, with $a,b\in A$ and $\omega\in A'$. It is clear that $V'\subseteq V$. 
\ssnl
On the other hand, take $a,b\in A$ and write
\begin{equation*}
\Delta(a)(1\ot b)=\sum_i p_i\ot q_i
\end{equation*}
with the $(q_i)$ linearly independent. Choose any index $j$ and a linear functional $\omega_j$ on $A$ that is $1$ and $q_j$ and $0$ on the other elements. Then
\begin{equation*}
p_j=(\iota\ot\omega_j)(\Delta(a)(1\ot b)).
\end{equation*}
We find that $p_j\in V'$. This holds for all indices and hence $\Delta(a)(1\ot b)\in V'\ot A$. Therefore $V'$ is the smallest subspace with this property for all $a,b$.
\ssnl
A similar argument works for the other cases.
\ebew

We relate these results with the notion of faithfulness of linear functionals.
\nl
\bf Faithful linear functionals  \rm 
\nl

\defin\label{defin:fff}
 A linear functional $\omega$ on $A$ is called {\it left faithful} if, given $a\in A$, we have $a=0$ if $\omega(ab)=0$ for all $b$ in $A$. Similarly, a linear functional on $A$ is called {\it right faithful} if, given $a\in A$, we have $a=0$ if $\omega(ca)=0$ for all $c\in A$. We call $\omega$ {\it faithful} if it is both left and right faithful.
\edefin 

In other words, $\omega$ is left faithful if $a\mapsto \omega(a\,\cdot\,)$ is an injective map from $A$ to $A'$ and it is right faithful if  $a\mapsto \omega(\,\cdot\,a)$ is injective.
\ssnl
Let us make a few remarks about this notion of faithfulness for linear functionals.

\opm
i) If the algebra is finite-dimensional, a linear functional that is left faithful will also be right faithful and vice versa. In other words, it is automatically faithful in that case. 
\ssnl
ii) For any faithful linear functional $\omega$ on a finite-dimensional algebra $A$, there exists an automorphism $\sigma$ of $A$ satisfying $\omega(ab)=\omega(b\sigma(a))$ for all $a,b\in A$. This automorphism is called the modular automorphism. Its inverse is the Nakayama automorphism, see \cite{Na}.
\ssnl
iii) In general, a faithful linear functional $\omega$ on an algebra that admits a modular automorphism as above is called weakly KMS (or simply a KMS-functional). Observe that it follows from the faithfulness that this modular automorphism is unique if it exists. But remark that there are examples of faithful linear functionals on infinite-dimensional non-degenerate algebras that do not have a modular automorphism, see e.g.\ \cite{VD-Ve}. 
\eopm

\oldcomment{Need some more comments on the terminology here? And references. Do we need these remarks about the modular automorphism? Do we want to refer to faithful sets, as formulated in Section 3?}{}

In Section \ref{s:ls-wmha} we will need the notion of a faithful set of linear functionals, see Definition \ref{defin:3.11c}.
\ssnl

As a consequence of Proposition \ref{prop:1.6c}  we now get the following. 

\prop\label{lem:span}
Assume that $\omega$ is  any linear functional on $A$. \
\ssnl
i) If $\omega$ is left faithful and if the map $T_1$ is regular, then the linear span of elements of the form $(\iota\ot \omega)(\Delta(a)(1\ot b))$, with $a,b\in A$, is equal to $V$. If it is right faithful and if $T_3$ is regular, the span of elements  $(\iota\ot \omega)((1\ot b)\Delta(a))$, with $a,b\in A$, is again equal to $V$. 
\ssnl
ii) If $\omega$ is right faithful and if $T_2$ is regular, the span of elements $(\omega\ot\iota)((c\ot  1)\Delta(a))$, with $a,c\in A$, is equal to $W$. If it is left faithful and if $T_4$ is regular, then the span of elements $(\omega\ot\iota)(\Delta(a)(c\ot  1))$, with $a,c\in A$, is also equal to $W$.
\eprop

\bew
Assume that $T_1$ is regular. By definition of $V$ we have $(\iota\ot\omega)(\Delta(p)(1\ot q))\in V$ for all $p,q\in A$ and all linear functionals $\omega\in A'$. 
Fix $\omega$ and denote by $V'$ the span of elements  $(\iota\ot\omega)(\Delta(p)(1\ot q))$ where $p,q\in A$. Then $V'\subseteq V$.
\ssnl
Now suppose that $\rho$ is an element in $A'$ such that $\rho$ is $0$ on $V'$. This means that $\omega(a)=0$ for all $a$ of the form $(\rho\ot\iota)(\Delta(p)(1\ot q))$. This holds for all $p,q$. Replace $q$ by $qc$. Then we get $\omega(ac)=0$ for all $c$ and all $a$ of the form $(\rho\ot\iota)(\Delta(p)(1\ot q))$. If now $\omega$ is left faithful, we must have $a=0$. Hence 
$(\rho\ot\iota)(\Delta(p)(1\ot q))=0$. This holds for all $p,q$ and therefore $\rho=0$ on $V$. This proves that $V\subseteq V'$. Because we already have the other inclusion, we get $V=V'$. 
\ssnl
This proves the first statement of item i). The other cases are proven in a similar way.
\ebew

Remark that we do not need that $\Delta$ is a homomorphism nor any form of coassociativity to define these legs of $\Delta$ and to obtain its properties.
\nl
\bf The notion of a counit \rm
\nl
Next we recall the notion of a counit in this setting. \oldcomment{\rood Again include a reference!}{}

\defin\label{defin:counit}
Let $\Delta$ be a linear map from $A$ to $M(A\ot A)$. Assume that $T_1$ and $T_2$ are regular. A linear functional $\varepsilon:A\to \mathbb C$ is called a \emph{counit} if
\begin{equation}(\varepsilon\ot\iota)(\Delta(a)(1\ot b))=ab
\tussenen
(\iota\ot\varepsilon)((c\ot 1)\Delta(a))=ca\label{eqn:counit}
\end{equation}
for all $a,b,c\in A$. 
\edefin

Again, we have the obvious definition in the other cases. And also here, the counit will be the same for the different possibilities.
\ssnl
Just as for the notion of fullness, also for the definition of a counit, there is no need to have that $\Delta$ is a homomorphism, nor that is coassociative. 

\opm
If the linear map $\Delta:A\to M(A\ot A)$ is full, then a counit is unique if it exists. In fact, in that case, if $\varepsilon$ is a linear map satisfying the first formula of (\ref{eqn:counit}) for all $a,b$ and if $\varepsilon'$ is a linear map satisfying the second formula of (\ref{eqn:counit}) for all $c,a$, then $\varepsilon=\varepsilon'$. This will follow from coassociativity. 
\ssnl
If the algebra $A$ has an identity, and if there is a counit, the map $\Delta$ is automatically full. In general, if there is a counit and if this counit is an algebra map, then again the coproduct is full. In particular, in these two cases, the counit is unique. In the general case however (e.g.\ for a weak multiplier Hopf algebra), there seems to be no way to conclude fullness from the existence of a counit.
\eopm

See Proposition 1.12 in \cite{VD-W0} for a proof of these statements.
\ssnl
Other notions needed, like an antipode and integrals, will be recalled later, where they are encountered in this note for the first time.

%%%%%%%%%%%%%%%%%%%%%%%%%%%%%%%%%%%%%%%%%%%%%%%%%%%%%%%%%%%%%%%%%%
%
% Bestand: artikel1.tex
%
%%%%%%%%%%%%%%%%%%%%%%%%%%%%%%%%%%%%%%%%%%%%%%%%%%%%%%%%%%%%%%%%%%

\section{\hspace{-17pt}. The Larson-Sweedler theorem for multiplier Hopf algebras}\label{s:ls-mha} %\input artikel2.tex

In this section, we start with a non-degenerate algebra $A$ and a coproduct $\Delta$ on $A$ as in Definition \ref{defin:copr} of Section \ref{s:prel}.  So it is a homomorphism $\Delta:A\to M(A\ot A)$ satisfying some form of coassociativity as in Definitions \ref{defin:coass} and \ref{defin:1.11a}. 
\ssnl
We introduce the notions of a left and a right integral. The main part of this section is devoted to the properties of the canonical maps $T_1$, $T_2$, $T_3$ and $T_4$ that can be proven from the existence of such integrals.
\nl
\bf Left and right Integrals\rm
\nl
In this section, we work with the following definition of a left and a right integral.

\defin\label{defin:int}
Consider a non-degenerate algebra $A$ with a coproduct $\Delta$ on $A$. 
\ssnl
i) A \emph{left integral} is defined if either $T_2$ or $T_4$ is regular as a non-zero linear functional $\varphi$ on $A$ satisfying
\begin{equation*}
(\iota\ot\varphi)((c\ot 1)\Delta(a)(c'\ot 1))=\varphi(a)cc'
\end{equation*}
for all $a,c,c'\in A$.
\ssnl
ii) A \emph{right integral} is defined if either $T_1$ or $T_3$ is regular as a non-zero linear functional $\psi$ on $A$ satisfying
\begin{equation*}
(\psi\ot \iota)((1\ot b')\Delta(a)(1\ot b))=\psi(a)b'b
\end{equation*}
for all $a,b,b'\in A$.
\edefin

Let us make some trivial remarks about these definitions. 

\opm\label{opm:2.2}
i) Let $\psi$ be a linear functional on $A$. Assume that $T_1$ is regular. Then $(\psi\ot\iota)\Delta(a)$ is defined as a left multiplier  on $A$ for all $a\in A$  by
\begin{equation*}
((\psi\ot\iota)\Delta(a))b=(\psi\ot\iota)(\Delta(a)(1\ot b))
\end{equation*}
where $b\in A$. When $\psi$ is non-zero it is a right integral if and only if $(\psi\ot\iota)\Delta(a)=\varphi(a)1$ for all $a$. If on the other hand $T_3$ is regular, then $(\psi\ot\iota)\Delta(a)$ is defined as a right multiplier  on $A$ for all $a\in A$  by
\begin{equation*}
b((\psi\ot\iota)\Delta(a))=(\psi\ot\iota)(1\ot b)\Delta(a))
\end{equation*}
where $b\in A$. When $\psi$ is non-zero, it is a right integral if and only if $(\psi\ot\iota)\Delta(a)=\varphi(a)1$ for all $a$.
In the event that both $T_1$ and $T_3$ are regular, then $(\psi\ot\iota)\Delta(a)$ is defined as a  multiplier of $A$. Again, if $\psi$ is non-zero, the multiplier is $\psi(a)1$ if and only if $\psi$ is a right integral.
\ssnl
ii) We have similar results for left integrals. If $T_2$ is regular, we have $(\iota\ot\varphi)\Delta(a)=\varphi(a)1$ as left multipliers while if $T_4$ is regular, we have this equation as right multipliers.
\eopm

All these properties  are direct consequences of the non-degeneracy of the product in $A$.
\ssnl
We agree that, when we assume the existence of a left integral, it is implicitly understood that either $T_2$ or $T_4$ is regular. Similarly, when we assume that a right integral exists, we implicitly require that $T_1$ or $T_3$  is regular.
\ssnl
Remark that it is not necessary that the linear map $\Delta:A\to M(A\ot A)$ is a homomorphism or that it is coassociative to define these integrals.
\snl
Before we continue, we want to add the following remark. It is important for the treatment of the Larson-Sweedler theorem for weak multiplier Hopf algebras in the next section, but it is a remark about integrals in the case of multiplier Hopf algebras.

\opm\label{opm:2.3}
Let $A$ be a \emph{unital} algebra and $\Delta:A\to A\ot A$ a coproduct on $A$.
Let $\varphi$ be a linear functional such that for all $a$ there is scalar $\lambda$ such that $(\iota\ot\varphi)\Delta(a)=\lambda 1$. 
Then we have 
\begin{align*}
\sum_{(a)} \Delta(a_{(1)})\,\varphi(a_{(2)})
&=\lambda\, 1\ot 1 \\
&=\sum_{(a)} (a_{(1)}\ot 1)\,\varphi(a_{(2)}).
\end{align*}
It follows that $(\iota\ot\varphi)\Delta(b)=\varphi(b)1$ for $b$ in the right leg of $\Delta$. Therefore, if $\Delta$ is full, we get it for all elements. Also if there is a counit, we can apply it from the very beginning and get  $\lambda=\varphi(a)$ so that also in this case, we get the formula for all elements.
\ssnl
An argument like this also works when $A$ is not unital but one has to be more  careful.  We have discussed this in an appendix because it is not relevant here. On the other hand, a similar problem occurs when treating the Larson-Sweedler theorem in the next section. There it is relevant however, see the item on integrals in Section \ref{s:ls-wmha}. 
\ssnl
In other words, we could give a more general definition of a left integral by requiring that $(\iota\ot\varphi)\Delta(a)$ is a scalar multiple of $1$ but we need that the right leg of $\Delta$ is all of $A$ to conclude that this is a left integral in the sense of Definition \ref{defin:int}.

\eopm
In what follows, \emph{we stick to the more common}, although somewhat less general, \emph{definition of integrals} as given in \ref{defin:int}.
\nl
To prove the main result, it will not be enough that non-zero integrals exist. We will need faithfulness of these integrals. We distinguish between left and right faithful as we defined in Definition \ref{defin:fff} in the previous section.
\nl
\bf Injectivity of the canonical maps \rm
\nl
First we prove results about the {\it injectivity of the canonical maps} when integrals exist. We distinguish between the different cases.

\prop\label{prop:Tinj}
i) The map $T_1$ is injective if there is a right faithful right integral.\\
ii) The map $T_2$ is injective if there is a left faithful left integral.\\
iii) The map $T_3$ is injective if there is a left faithful  right integral.\\
iv) The map $T_4$ is injective if there is a right faithful left integral.
\eprop

\bew
i) Take an element $z\in A\ot A$ and assume that $T_1(z)=0$. Write $z=\sum_i a_i\ot b_i$ so that 
\begin{equation}
\sum_i \Delta(a_i)(1\ot b_i)=0. \label{eqn:Tz0}
\end{equation}
Take any $x,y\in A$, multiply  with $(1\ot y)\Delta(x)$ from the left and apply a right integral $\psi$ on the first leg. Using that $\Delta$ is a homomorphism, we get 
\begin{equation*}
\sum_i \psi(xa_i) yb_i=0.
\end{equation*}
This holds for all $x$ and if $\psi $ is right faithful we get $\sum_i a_i\ot yb_i=0$. This holds for all $y$ and so $\sum_i a_i\ot b_i=0$.  Hence $T_1$ is injective.
\ssnl
We have implicitly used that $T_1$ or $T_3$ is regular by the existence of $\psi$.
\ssnl
ii) Next assume $\sum_i (c_i\ot 1)\Delta(a_i)=0$. Multiply with $\Delta(x)(y\ot 1)$ from the right and apply a left integral $\varphi$. We get
 $\sum_i \varphi(a_ix)c_iy=0$. If $\varphi$ is left faithful, we find $\sum_i c_i\ot a_i=0$ and so $T_2$ is injective.
 \ssnl
Properties  iii) and iv) are proven in the same way.
\ebew

The proof is simple and is already found in \cite{VD-W-ls}. We recall it here to verify the result when we use only one-sided faithfulness of the integrals. Also observe that for this result, regularity of the coproduct is not needed, also coassociativity is not needed, but we need $\Delta$ to be a homomorphism.
\nl
\bf Surjectivity of the canonical maps \rm
\nl
We will now further concentrate on the {\it ranges of the canonical maps}. As for the results on injectivity, also here we carefully distinguish between the different cases. 
\ssnl
Before we start, let us consider the case of a unital algebra $A$. We will use it to motivate the choices made further.
\voorb
Let $A$ be a unital algebra and $\Delta:A\to A\ot A$ a coproduct on $A$. Assume that $\varphi$ is a left integral. Take  elements $p,q,b\in A$ and consider 
\begin{equation}
y=\sum_{(p),(q)}(p_{(1)}\ot q_{(1)}b)\,\varphi(p_{(2)}q_{(2)}).\label{eqn:2.2}
\end{equation}
Using coassociativity, that $\Delta$ is a homomorphism and $\varphi$ left invariant, we have
\begin{equation*}
T_1(y)
=\sum_{(p),(q)}(p_{(1)}\ot p_{(2)}q_{(1)}b)\,\varphi(p_{(3)}q_{(2)})
=\sum_{(p)}(p_{(1)}\ot b)\,\varphi(p_{(2)}q).
\end{equation*}
So $T_1(y)=a\ot b$ where $a=(\iota\ot\varphi)(\Delta(p)(1\ot q))$. 
\ssnl
If $\varphi$ is faithful, the span of such elements $a$ is equal to the left leg of $\Delta$. In particular, if $\Delta$ is full, we see that the range of $T_1$ is all of $A\ot A$.
\evoorb

One can obtain similar results for the other canonical maps. And it is not hard to obtain the Larson-Sweedler result for unital algebras from this.
\ssnl
In fact, to obtain the result in the case of multiplier Hopf algebras (as in \cite{VD-W-ls}) and in the case of weak multiplier Hopf algebras (as in \cite{K-VD}), the same kind of elements are used. In these papers, the coproduct is assumed to be regular so that all canonical maps are regular. Then one can use the formula (\ref{eqn:2.2}) to define an element $y$ in $A\ot A$ because $q_{(1)}$ is covered by $b$ and consequently $p_{(2)}$ is also covered.
\ssnl
In what follows we now modify the arguments at a few points. In the first place, we do not assume that all canonical maps are regular, use one-sided faithfulness of the integrals where possible and finally, we do not require the coproduct to be full. 
\ssnl
We will systematically treat the four canonical maps. As before, we assume that $A$ is a non-degenerate algebra and $\Delta$ a coproduct on $A$. We consider a left integral $\varphi$ and a right integral $\psi$.

\prop\label{prop:rTl} 
i) Assume that the canonical maps $T_1$ and $T_4$ are regular.  Given $p,q,b\in A$ we can define an element  $y\in A\ot A$ by
\begin{equation*}
y=(\iota\ot\iota\ot\varphi)(\Delta_{13}(p)\Delta_{23}(q)(1\ot b\ot 1)).
\end{equation*}	 
We can apply $T_1$ and we obtain $T_1(y)=a\ot b$ where $a=(\iota\ot \varphi)(\Delta(p)(1\ot q))$. 
\ssnl
ii) Assume that the canonical maps $T_2$ and $T_3$ are regular. Given $p,q,c\in A$ we can define an element $y\in A\ot A$ by
\begin{equation*}
y=(\psi\ot\iota\ot\iota)((1\ot c\ot 1)\Delta_{12}(p)\Delta_{13}(q)).
\end{equation*}
We can apply $T_2$ and we obtain $T_2(y)=c\ot a$ where  $a=(\psi\ot\iota)((p\ot 1)\Delta(q))$. 
\ssnl
iii) Assume that the canonical maps $T_2$ and $T_3$ are regular. Given $p,q,b\in A$ we can define an element $z\in A\ot A$ by
\begin{equation*}
z=(\iota\ot\iota\ot\varphi)((1\ot b\ot 1)\Delta_{23}(p)\Delta_{13}(q)).
\end{equation*}
We can apply $T_3$ and we obtain $T_3(z)=a\ot b$ where now $a=(\iota\ot\varphi)((1\ot p)\Delta(q))$.
\ssnl
iv) Finally assume that the maps $T_1$ and $T_4$ are regular. Given $p,q,c\in A$ we can define $z\in A\ot A$ by
\begin{equation*}
z=(\psi\ot\iota\ot\iota)(\Delta_{13}(p)\Delta_{12}(q)(1\ot c\ot 1)).
\end{equation*}
We can apply $T_4$ and we obtain $T_4(z)=c\ot a$ where $a=(\psi\ot\iota)(\Delta(p)(q\ot 1))$. 
\eprop

\bew
i) Assume  that $T_1$ and $T_4$ are regular. Take $p,q,b\in A$ and define $U$ in the multiplier algebra $M(A\ot A\ot A)$ by
\begin{equation*}
U=\Delta_{13}(p)\Delta_{23}(q)(1\ot b\ot 1).
\end{equation*}
Because $T_4$ is regular,  $\Delta(q)(b\ot 1)$ belongs to $A\ot A$ and we can write it as $\sum_i r_i\ot s_i$. Then we find
\begin{equation*}
U=\sum_i \Delta_{13}(p)(1\ot r_i\ot s_i).
\end{equation*}
Because $T_1$ is regular, $\Delta(p)(1\ot s_i)$ belongs to $A\ot A$ for all $i$ and we find that $U\in A\ot A\ot A$. We apply $T_1\ot \iota$ and we get 
\begin{equation*}
(T_1\ot\iota)(U)=\sum_i (\Delta\ot\iota)(\Delta(p)(1\ot s_i))(1\ot r_i\ot 1).
\end{equation*}
We multiply with an element $c$ of $A$ from the right  in the first factor. Then we can use coassociativity for the combination of the canonical maps $T_1$ and $T_4$ (see Definition \ref{defin:1.11a}), and we obtain
\begin{align*}
((T_1\ot\iota)(U))(c\ot 1\ot 1)
&=\sum_i (\Delta\ot\iota)(\Delta(p)(1\ot s_i))(c\ot 1\ot 1)(1\ot r_i\ot 1)\\
&=\sum_i (\iota\ot\Delta)(\Delta(p)(c\ot 1))(1\ot 1\ot  s_i)(1\ot r_i\ot 1)\\
&=(\iota\ot\Delta)(\Delta(p)(c\ot 1))(1\ot \Delta(q))(1\ot b\ot 1)\\
&=(\iota\ot\Delta)(\Delta(p)(c\ot 1)(1\ot q))(1\ot b\ot 1)\\
&=(\iota\ot\Delta)(\Delta(p)(1\ot q))(c\ot b\ot 1).
\end{align*}
We have used that $\Delta$ is a homomorphism. In the end we can cancel $c$ to obtain
\begin{equation*}
(T_1\ot\iota)(U)=(\iota\ot\Delta)(\Delta(p)(1\ot q))(1\ot b\ot 1).
\end{equation*}
Finally, we apply the left integral on the third factor and we get, because $(\iota\ot\varphi)U=y$, and using left invariance of $\varphi$ that
\begin{equation*}
T_1(y)=(\iota\ot\varphi)(\Delta(p)(1\ot q))\ot b.
\end{equation*}
So $T_1(y)=a\ot b$ when we define $a=(\iota\ot\varphi)(\Delta(p)(1\ot q))$. This completes the proof of the first statement in the proposition.
\ssnl
ii) The proof is completely similar to the proof of the previous result.  But let us give a proof using the Sweedler notation, just to illustrate how that works.
\ssnl
The element $y$ is given by
\begin{equation*}
y=\sum_{(p)(q)} \psi(p_{(1)}q_{(1)})\, cp_{(2)}\ot q_{(2)}.
\end{equation*}
The element $c$ covers $p_{(2)}$ (using that $T_3$ is regular) and then $p_{(1)}$ covers $q_{(1)}$ (using that $T_2$ is regular).  We can apply the map $T_2$ to get
\begin{equation*}
T_2(y)= \sum_{(p)(q)} \psi(p_{(1)}q_{(1)})\, cp_{(2)}q_{(2)}\ot q_{(3)}.
\end{equation*}
Again $c$ will cover $p_{(2)}$ (as $T_3$ is regular)  and consequently also $q_{(1)}$ and $q_{(2)}$ are covered (as $T_2$ is regular). This requires another form of coassociativity however. In order to use coassociativity with the canonical maps $T_2$ and $T_3$ as in Definition \ref{defin:1.11a} we need an extra covering of  $q_{(3)}$ in the above formula with an element $b$ on the left. We also had to  do this in the previous proof.  
Finally, we can use right invariance of $\psi$ to get
\begin{equation*}
T_2(y)= \sum_{(q)} \psi(pq_{(1)})\, c\ot q_{(2)}
\end{equation*}
and we get $T_2(y)=c\ot a$ with $a= \sum_{(q)} \psi(pq_{(1)})q_{(2)}$.
\ssnl
iii)  If now $T_2$ and $T_3$ are regular, we can define the element $z\in A\ot A$. A similar calculation as above will give
\begin{equation*}
T_3(z)=(\iota\ot\varphi)((1\ot p)\Delta(q))\ot b.
\end{equation*}
\ssnl
iv) Again if $T_1$ and $T_4$ are regular, we can define the element $z\in A\ot A$ as in the formulation and obtain the desired formula.
\ebew

The reader is advised to compare the proof of item ii), using the Sweedler notation, with the proof of item i) where we did not use the Sweedler notation.
\ssnl
We used that $\Delta$ is a homomorphism and that it is coassociative. However, we do not need any faithfulness property of the integrals for these  results. Of course, without this condition on the integrals, we have no guarantee that the range of these canonical maps is big enough.
\ssnl
If we assume that the integrals are faithful and that the coproduct is full, we get from the above results that the canonical maps are surjective. However, we do not assume this here and so we have to be more careful. 
\nl
\bf The main result \rm
\nl
We now combine the results of Proposition \ref{prop:Tinj} with those of Proposition \ref{prop:rTl} to obtain the main results for multiplier Hopf algebras.
\ssnl
Recall from Definition \ref{defin:full-cp} that the left and the right leg of $\Delta$ are the smallest subspaces $V$ and $W$ of $A$ satisfying
\begin{align*}
&\Delta(a)(1\ot b)\in V\ot A
\ \tussenen 
(1\ot b)\Delta(a)\in V\ot A\\
&\Delta(a)(c\ot 1)\in A\ot W
\tussenen
(c\ot 1)\Delta(a)\in A\ot W
\end{align*}
for all $a,b,c\in A$. It is implicitly assumed that the canonical maps, involved in these formulas, are regular.
\ssnl
We do not assume that $\Delta$ is full. So the spaces $V$ and $W$ could be strictly smaller than $A$. However, as we mentioned already, it will eventually be shown that they are equal to $A$. This will follow from the other assumptions.
\oldcomment{We need the two integrals to obtain these results. In the conclusion sections, we need to refer to this as Kusterman's trick.}{}

\prop\label{prop:2.7c}
i) If $T_1$ and $T_4$ are regular, then $T_1$ is bijective if there is a left faithful left integral and a right faithful right integral. \\
ii) If $T_2$ and $T_3$ are regular, then $T_2$ is bijective if there is a left faithful left integral and a right faithful right integral.\\
iii) If $T_2$ and $T_3$ are regular, then $T_3$ is bijective if there is a right faithful left integral and a left faithful right integral.\\
iv) If $T_1$ and $T_4$ are regular, then $T_4$ is bijective if there is a right faithful left integral and a left faithful right integral.
\eprop

\bew
i) Because we have a right faithful right integral, by Proposition \ref{prop:Tinj} the map $T_1$ is injective. To prove that it is surjective, we proceed as follows.
\ssnl
Take any element $x\in A\ot A$. We know that $T_1(x)\in V\ot A$. On the other hand, because $\varphi$ is assumed to be left faithful, by Proposition  \ref{lem:span}, $V$ is spanned by elements of the form $a=(\iota\ot\varphi)((\Delta(p)(1\ot q))$ for $p,q\in A$. In Proposition \ref{prop:rTl} we have seen that for such elements we have 
$a\ot b=T_1(y)$ where 
\begin{equation*}
y= \sum_{(p),(q)}p_{(1)}\ot q_{(1)}b\, \varphi(p_{(2)}q_{(2)}).
\end{equation*}
The element $y$  is a linear combination of elements of the form $\sum_{(p)}p_{(1)}\ot r\, \varphi(p_{(2)}s)$. Hence $y\in V\ot A$.
Taking all these results together, we find that $T_1(x)\in T_1(V\ot A)$. Because $T_1$ is injective, it follows that $x\in V\ot A$. This is now true for all elements $x\in A\ot A$ and hence $V=A$. Then also $T_1$ is surjective.
\ssnl
ii) The other statements are proven in a similar way by using the other parts of Propositions \ref{prop:Tinj}  and \ref{prop:rTl}.
\ebew

Observe that replacing $\Delta$ by $\Delta^\text{cop}$ results in interchanging $T_1$ with $T_4$ and a left integral with a right integral. This shows that i) and iv) are compatible with each other. It will also interchange $T_2$ with $T_3$ and so also ii) and iii) are compatible with each other. On the other hand, if we replace $A$ by $A^\text{op}$, this converts $T_1$ to $T_3$ and $T_2$ to $T_4$, while left faithfulness becomes right faithfulness. This shows that i) is compatible with iii) and that ii) is compatible with iv).
\oldcomment{Dit is terug helemaal nagekeken.}{}
\ssnl
Also remark that  in all these cases, we need regularity of a second canonical map to obtain a result about one canonical map. We see e.g. that we need also regularity of $T_4$ in order to prove that $T_1$ is bijective. Similarly for the other three cases.
\ssnl
The argument we use here to obtain the surjectivity of the canonical maps is inspired by an argument given in \cite{Ku-Va}, known as `Kustermans trick', see a comment in Section \ref{s:concl} and further Appendix \ref{sect:appB} where we say something more about this.
 
 \opm\label{opm:2.8d}
 i) We see that we get a full coproduct under the appropriate conditions. \\
 %\ssnl
 ii) We can also obtain that $A$ is idempotent if e.g.\ $T_1$ and $T_4$ are regular and if any of these canonical maps is surjective. Indeed, assume that $\omega=0$ on $A^2$. If $T_4$ is regular, we have $$\Delta(a)(p\ot q)\in A\ot A^2$$ and so $(\iota\ot\omega)(\Delta(a)(p\ot q))=0$ for all $a,p,q$. If also $T_1$ is regular we can conclude that still $(\iota\ot\omega)(\Delta(a)(1\ot q))=0$ for all $a,q$. So when $T_1$ is surjective, we must have $\omega=0$.
 This proves that $A=A^2$.\\
 iii) Assume that $A$ is idempotent. If $T_1(A\ot A)=A\ot A$ we see that also $\Delta(A)(A\ot A)=A\ot A$. This also holds if $T_4$ is surjective. On the other hand, if $T_2$ is surjective we would get $(A\ot A)\Delta(A)=A\ot A$. In these cases we get that $\Delta$ is a non-degenerate homomorphism from $A$ to $M(A\ot A)$.
 \eopm
 
 We now come to our main results. 
 
 \stel \label{stel:2.8}
 Assume that $A$ is a non-degenerate algebra and that $\Delta:A\to M(A\ot A)$ is a \emph{regular} coproduct on $A$. Then $(A,\Delta)$ is a multiplier Hopf algebra if there is a left faithful left integral and a right faithful right integral. It is a regular multiplier Hopf algebra if there is a faithful left integral and a faithful right integral.
 \estel
 
\bew
i) Because we assume that $\Delta$ is regular, all canonical maps are regular. And because we assume that there is a left faithful left integral and a right faithful right integral, we can apply items i) and ii) of the previous proposition and we have that $T_1$ and $T_2$ are bijective maps. Hence, by definition, $(A,\Delta)$ is a multiplier Hopf algebra.
\ssnl
ii) If there is a faithful left integral and a faithful right integral, also the other items, iii) and iv) apply and the canonical maps $T_3$ and $T_4$ are also bijective. This means that $(A,\Delta)$ is a regular multiplier Hopf algebra.
\ebew
 
The last property is obtained in  \cite{VD-W-ls} under the extra assumption that the coproduct is full. We see that this result is a bit stronger than the one in \cite{VD-W-ls}. Observe that the coproduct in any multiplier Hopf algebra is automatically full.
 \ssnl
 In the following result we get a still stronger property than the last one.   For the notion of a single sided (left or right) multiplier Hopf algebra, we refer to \cite{VD-ssmhas}. 
 
 \prop\label{prop:2.9d} 
 Assume that $A$ is a non-degenerate algebra and that $\Delta$ is a coproduct on $A$. If $T_1$ and $T_4$ are regular, then $(A,\Delta)$ is a left multiplier Hopf algebra if there is a faithful left integral and a faithful right integral. If $T_2$ and $T_3$ are regular, then $(A,\Delta)$ is a right multiplier Hopf algebra if there is a faithful left integral and a faithful right integral.
 \eprop
 
 \bew
 i) Given that $T_1$ and $T_4$ are regular, we can use item i) and item iv) of Proposition \ref{prop:2.7c}. The maps $T_1$ and $T_4$ are bijective because there is a faithful left integral and a faithful right integral. Then $(A,\Delta)$ is a left multiplier Hopf algebra as in item i) of Definition 1.5 of \cite{VD-ssmhas}.
 \ssnl
 ii) Given that $T_2$ and $T_3$ are regular, we can use item ii) and item iii) of Proposition \ref{prop:2.7c}. Again the maps $T_2$ and $T_3$ are bijective because of the existence of faithful integrals. Now $(A,\Delta)$ is a right multiplier Hopf algebra as in item ii) of Definition 1.5 of \cite{VD-ssmhas}.
 \ebew

 For the above result we only need the material obtained in this paper. However, it is shown in \cite{VD-ssmhas}  that a single sided multiplier Hopf algebra is automatically a regular multiplier Hopf algebra. Therefore, we also get the following stronger form of the second part of Theorem \ref{stel:2.8}.
 
\stel\label{stel:2.11} 
Assume that $A$ is a non-degenerate algebra and that $\Delta:A\to M(A\ot A)$ is a coproduct on $A$. If $T_1$ and $T_4$ are regular, then $(A,\Delta)$ is a regular multiplier Hopf algebra if there is a faithful left integral and a faithful right integral. If $T_2$ and $T_3$ are regular, then $(A,\Delta)$ is a regular multiplier Hopf algebra if there is a faithful left integral and a faithful right integral.
 \estel
 
 We now look at some consequences for the case when $A$ is unital. We derive the following consequence from Theorem \ref{stel:2.8}.
 
 \gev
 Assume that $A$ is a \emph{unital} algebra and that $\Delta:A\to A\ot A$ is a coproduct on $A$. Then $A$ is a Hopf algebra if there exists a left faithful left integral and a right faithful right integral. If these integrals are faithful, $A$ has an invertible antipode.
 \egev
 
 \bew
Because $A$ has an identity all the canonical maps are automatically regular. 
\ssnl
i) If there is a left faithful left integral and a right faithful right integral, the maps $T_1$ and $T_2$ are bijections. Then $A$ is a Hopf algebra.
\ssnl
ii) If there is a faithful left integral and a faithful right integral, all four canonical maps are bijective. Then $A$ is a Hopf algebra with an invertible antipode.
\ebew

Observe that we do not need the existence of a counit to obtain this result. It is a consequence of the other conditions.
\ssnl
 For a finite-dimensional space, a left faithful functional is automatically faithful. Similarly for a right faithful functional. Moreover a finite-dimensional algebra with a faithful functional is automatically unital. Hence we obtain the following consequence.
 
 \gev
 Assume that $A$ is a finite-dimensional algebra and that  $\Delta:A\to A\ot A$ is a coproduct on $A$. Then $A$ is a Hopf algebra if there exists a  faithful left integral and a  faithful right integral. 
 \egev
 
 We also obtain that the antipode is invertible, but that is always the case for a finite-dimensional Hopf algebra.
\nl
\bf Further reflections \rm
\nl
From the general theory, we know that for any multiplier Hopf algebra there exists a counit and an antipode. The result is also true for single sided multiplier Hopf algebras. In the first case, this is already proven in \cite{VD-mha}. In the second case, a somewhat easier approach is found in \cite{VD-ssmhas}. In the two cases, these objects are constructed using the bijectivity of the canonical maps.
\ssnl
In this paper, the bijectivity of the canonical maps is proven  by using properties of the integrals. In this item, we show how the counit and antipode can be constructed directly from the properties of the integrals, thus avoiding the intermediate step. It turns out to give also some extra information.
\nl
We begin with the counit. In \cite{VD-ssmhas} the counit is obtained from the inverse of the canonical map $T_1$. If $a,b\in A$ and 
\begin{equation*}
a\ot b=\sum_i \Delta(p_i)(1\ot q_i),
\end{equation*}
then $\varepsilon(a)b=\sum_i p_iq_i$. This however only defines $\varepsilon(a)$ as a left multiplier of $A$. The canonical map $T_4$ is used to show that this left multiplier is actually a scalar multiple of the identity.
\ssnl
Here we proceed as follows.

\prop
Let $A$ be a non-degenerate algebra and $\Delta:A\to M(A\ot A)$ a coproduct on $A$. Assume that the maps $T_1$ and $T_4$ are regular. Also assume that we have a left faithful left integral and a right faithful right integral. Then one can define a linear map $\varepsilon$ on the left leg of $\Delta$ satisfying 
$$(\varepsilon\ot\iota)(\Delta(a)(1\ot b))=\varepsilon(a)b$$
 for all $a,b$.
\eprop
\bew
i) Because $\varphi$ is assumed to be left faithful, the left leg $V$ of $\Delta$ is spanned by elements of the form $(\iota\ot\varphi)(\Delta(p)(1\ot q))$, see item i) in Proposition \ref{lem:span}.
\ssnl
ii) We want to define $\varepsilon$ on $V$ by the formula
\begin{equation*}
\sum_i \varepsilon((\iota\ot\varphi)(\Delta(p)(1\ot q)))=\varphi(pq)
\end{equation*}
where $p,q\in A$. To show that this is  well-defined, we first assume that we have elements $p_i, q_i$ so that 
\begin{equation*}
\sum_i (\iota\ot\varphi)(\Delta(p_i)(1\ot q_i))=0.
\end{equation*}
By the first result in Proposition \ref{prop:rTl} it follows that 
\begin{equation*}
T_1(\sum_i (\iota\ot\iota\ot\varphi)(\Delta_{13}(p_i)\Delta_{23}(q_i)(1\ot b\ot 1)))=0
\end{equation*}
for all $b$. Because we have a right faithful right integral, we can apply  Proposition \ref{prop:Tinj} and the canonical map $T_1$ is injective. So
\begin{equation*}
\sum_i (\iota\ot\iota\ot\varphi)(\Delta_{13}(p_i)\Delta_{23}(q_i)(1\ot b\ot 1))=0
\end{equation*}
for all $b$. If we now apply multiplication and use left invariance of $\varphi$, we find
\begin{equation*}
0=\sum_i (\iota\ot\varphi)(\Delta(p_i)\Delta(q_i)(b\ot 1))=\sum_i b\,\varphi(p_iq_i)
\end{equation*}
for all $b$ and so $\sum_i \varphi(p_iq_i)=0$. This means that we can define $\varepsilon$ on $A$ by the formula
\begin{equation*}
\varepsilon((\iota\ot\varphi)(\Delta(p)(1\ot q)))=\varphi(pq).
\end{equation*}
\ssnl
iii) For any $a,b,c\in A$ we have
\begin{equation*}
\varepsilon((\iota\ot\varphi)(\Delta(a)(1\ot bc)))=\varphi(abc)
\end{equation*}
and because $\varphi$ is assumed to be left faithful, it follows that 
\begin{equation*}
(\varepsilon\ot\iota)((\Delta(a)(1\ot b))=ab.
\end{equation*}

\vskip-12pt
\ebew

The proof uses that $T_4$ is regular by the use of Proposition \ref{prop:rTl}. In the existence proof in \cite{VD-ssmhas}, it is needed that it is surjective, see Proposition 2.2 in that paper. On the other hand, here we only have the counit on the left leg of $\Delta$. This is not really a problem because, after all, the formula characterizing the counit only uses the values on the left leg of $\Delta$.
Apart from this difference, there is also some similarity. In \cite{VD-ssmhas} the counit is defined by $\varepsilon(a)b=mT_1\inv(a\ot b)$. Also here we use in fact the same formula. 
\ssnl
To continue, we use the same method to obtain a counit $\varepsilon'$ on the right leg of $\Delta$.

\prop\label{prop:2.14}
Let $A$ be a non-degenerate algebra and $\Delta:A\to M(A\ot A)$ a coproduct on $A$. Assume that the maps $T_1$ and $T_4$ are regular. Also assume that we have a right faithful left integral and a left faithful right integral. Then one can define a linear map $\varepsilon'$ on the right leg of $\Delta$ satisfying 
$$(\iota\ot\varepsilon')(\Delta(a)(c\ot 1))=\varepsilon(a)c$$
 for all $a,c$.
\eprop
\bew
We define $\varepsilon'$ by
\begin{equation*}
\varepsilon'((\psi\ot\iota)(\Delta(p)(q\ot 1))=\psi(pq)
\end{equation*}
and we use a similar argument as in the previous proof to show that it is well-defined. We need regularity of $T_1$ and $T_4$ to use Proposition \ref{prop:rTl} and a right faithful left integral to use that $T_4$ is injective. Then we have 
\begin{equation*}
\varepsilon'((\psi\ot\iota)(\Delta(p)(qx\ot 1))=\psi(pqx)
\end{equation*}
and using the left faithfulness of $\psi$ that 
\begin{equation*}
(\iota\ot\varepsilon')(\Delta(p)(q\ot 1))=pq.
\end{equation*}

\vskip -12pt \ebew

The proof uses that $T_1$ is regular. In the existence proof in \cite{VD-ssmhas}, it is needed that $T_1$ is surjective. But just as in that paper, we essentially define here $\varepsilon'$ by $$\varepsilon'(a)c=m^\text{op} T_4\inv(c\ot a).$$
\ssnl
Observe that, in order to get $\varepsilon$ and $\varepsilon'$ we need a faithful left integral and a faithful right integral. This is in agreement with the result of Proposition \ref{prop:2.9d}.
\ssnl
The case of the counit gives little extra information, but we include it to illustrate the approach.
\nl
We now consider the antipode.
\ssnl
First recall the following well-known result.

\prop\label{prop:Sform}
Let $(A,\Delta)$ be a regular multiplier Hopf algebra with antipode $S$. If $\varphi$ is a left integral, we have 
\begin{equation}
S((\iota\ot\varphi)(\Delta(a)(1\ot b)))=(\iota\ot\varphi)((1\ot a)\Delta(b))\label{eqn:Sform}
\end{equation}
for all $a,b\in A$.
\eprop

With the use of the Sweedler notation, we get
\begin{equation*}
(1\ot a)\Delta(b)=\sum_{(a)}(S(a_{(1)})\ot 1)\Delta(a_{(2)}b)
\end{equation*}
and the formula (\ref{eqn:Sform}) follows by applying $\varphi$ on the second leg of this equation.
\ssnl
Consider a faithful left integral $\varphi$. Then the antipode $S$, if it exists, is completely determined by the formula (\ref{eqn:Sform}) of Proposition \ref{prop:Sform}. We now use this formula to define the antipode in this context.

\prop\label{prop:2.26}
Let $A$ be a non-degenerate algebra and $\Delta:A\to M(A\ot A)$ a coproduct. Assume that the maps $T_1$ and $T_4$ are regular. Also assume that we have a left faithful left integral and a right faithful right integral. Then we can define a linear map $S$ from $A$ to the left multiplier algebra $L(A)$ satisfying  
\begin{equation}
S((\iota\ot\varphi)(\Delta(p)(1\ot q)))=(\iota\ot\varphi)((1\ot p)\Delta(q)).\label{eqn:2.12}
\end{equation}
as left multipliers.
\eprop
\bew
We proceed as in the proof of Proposition \ref{prop:2.14}. 
\ssnl
i) Because $\varphi$ is assumed to be left faithful, the left leg $V$ of $\Delta$ is spanned by elements of the form $(\iota\ot\varphi)(\Delta(p)(1\ot q))$, see item i) in Proposition \ref{lem:span}.
\ssnl
ii) We now want to define $S$ by the formula
\begin{equation*}
S((\iota\ot\varphi)(\Delta(p)(1\ot q)))=(\iota\ot\varphi)((1\ot p)\Delta(q)).
\end{equation*}
To do this, assume that $(p_i)$ and $(q_i)$ are elements in $A$ so that\\ $\sum_i (\iota\ot\varphi)(\Delta(p_i)(1\ot q_i))=0$. Then by Proposition \ref{prop:rTl} we have
\begin{equation*}
T_1(\sum_i (\iota\ot\iota\ot\varphi)(\Delta_{13}(p_i)\Delta_{23}(q_i)(1\ot b\ot 1)))=0
\end{equation*}
for all $b$ in $A$. Because we assume that there is a right faithful right integral, by the injectivity of the canonical map $T_1$  we find
\begin{equation*}
\sum_i (\iota\ot\iota\ot\varphi)(\Delta_{13}(p_i)\Delta_{23}(q_i)(1\ot b\ot 1))=0
\end{equation*}
for all $b$. We can apply the counit on the first factor and we obtain
\begin{equation*}
\sum_i (\iota\ot\varphi)((1\ot p_i)\Delta(q_i)(b\ot 1))=0.
\end{equation*}
We see that we can define $S(a)$ for $a\in A$ as a left multiplier by the formula
\begin{equation*}
S((\iota\ot\varphi)(\Delta(p)(1\ot q)))=(\iota\ot\varphi)((1\ot p)\Delta(q)).
\end{equation*}

\vskip -12pt
\ebew

\prop Assume the conditions of the previous proposition. 
For all $a$ in the left leg of $\Delta$ and for all $b$ we have $\sum_{(a)} a_{(1)}\ot S(a_{(2)})b$ well-defined in $A\ot A$ and $\sum_{(a)} a_{(1)}S(a_{(2)})b=\varepsilon(a)b$.
\eprop
\bew
Apply the formula (\ref{eqn:2.12})with $p$ replaced by $(\omega\ot\iota)(\Delta(p)(c\ot 1))$ where $p,b,c\in A$ and $\omega\in A'$. Using the Sweedler notation we find
\begin{equation*}
\sum_{(p)} \omega(p_{(1)}c)S(p_{(2)})b\,\varphi(p_{(3)}q)
 =\sum_{(p),(q)} \omega(p_{(1)}c)q_{(1)}b\,\varphi(p_{(2)}q_{(2)}).
\end{equation*}
This means that 
\begin{equation*}
\sum_{(p)} p_{(1)}c\ot S(p_{(2)})b\,\varphi(p_{(3)}q)
 =\sum_{(p),(q)}p_{(1)}c\ot q_{(1)}b\,\varphi(p_{(2)}q_{(2)}).
\end{equation*}
We can write this as
\begin{equation}
\sum_{(a)} a_{(1)}c\ot S(a_{(2)})b
 =\sum_{(p),(q)}p_{(1)}c\ot q_{(1)}b\,\varphi(p_{(2)}q_{(2)})\label{eqn:2.5d}
\end{equation}
where $a=(\iota\ot\varphi)(\Delta(p)(1\ot q))$. 
On the right hand side, we see that $b$ covers $q_{(1)}$ and also $p_{(2)}$ is covered. We can cancel $c$ and still get an element in $A\ot A$. It implies that $\sum_{(a)} a_{(1)}c\ot S(a_{(2)})b$ is well-defined in $A\ot A$. Then we can apply multiplication and obtain
\begin{align*}
\sum_{(a)} a_{(1)}S(a_{(2)})b)
& =\sum_{(p),(q)}p_{(1)}q_{(1)}b\varphi(p_{(2)}q_{(2)})=b\varphi(pq)\\
& =b\varepsilon((\iota\ot\varphi)(\Delta(p)(1\ot q)).
\end{align*}
We see that $\sum_{(a)}a_{(1)}S(a_{(2)})b=\varepsilon(a)b$. 
\ebew

The map $a\ot b\mapsto \sum_{(a)}a_{(1)}\ot S(a_{(2)})b$ is the inverse of $T_1$. One can see this by applying $T_1$ on the formula (\ref{eqn:2.5d}) above.
\nl

We finish with some remarks about the more general situations. 
\ssnl
In Theorem \ref{stel:2.8} we found a multiplier Hopf algebra under the conditions
(1) that the coproduct is regular, (2) the existence of a left faithful left integral and (3) the existence of a right faithful right integral.
We see from the formula (\ref{eqn:Sform}), that indeed the antipode will map $A$ into $A$ when it is assumed that the coproduct is regular. Now, for a general multiplier Hopf algebra, the antipode is known to map into the multiplier algebra $M(A)$. For this to happen, we should weaken the condition on the coproduct and only require the canonical maps $T_1$ and $T_2$ to be regular. 

%%%%%%%%%%%%%%%%%%%%%%%%%%%%%%%%%%%%%%%%%%%%%%%%%%%%%%%%%%%%%%%%%%%%%%%
%
% Bestand: artikel3.tex
% 
% %%%%%%%%%%%%%%%%%%%%%%%%%%%%%%

\section{\hspace{-17pt}. The theorem for  weak multiplier Hopf algebras} \label{s:ls-wmha} % \input artikel3.tex

In this section, we consider the Larson-Sweedler theorem for {\it weak} multiplier Hopf algebras.

\ssnl
We formulate  the starting conditions, compare them with the ones we had in the previous section, as well as with the conditions we used in the original article \cite{K-VD}. It is the intention to weaken these conditions where possible.
\nl
\bf Conditions on the algebra $A$ and the coproduct $\Delta$ \rm
\nl
We assume that $A$ is a non-degenerate algebra and that $\Delta:A\to M(A\ot A)$ is a coproduct (as in Definition \ref{defin:copr}).
So $\Delta$ is a homomorphism and it is coassociative in any of the forms we have in  Definition \ref{defin:coass} and Definition \ref{defin:1.11a}. 
\ssnl
Just as in \cite{K-VD} we have the following condition on the coproduct.

\voorw\label{voorw:3.1c}
We assume the existence of an idempotent $E\in M(A\ot A)$ satisfying
\begin{equation*}
\Delta(A)(A\ot A)=E(A\ot A) 
\tussenen
(A\ot A)\Delta(A)=(A\ot A)E.
\end{equation*}
\evoorw

We know that such an idempotent is unique and it is the smallest one satisfying
\begin{equation*}
E\Delta(a)=\Delta(a)E=\Delta(a)
\end{equation*}
for all $a\in A$.
One can extend $\Delta\ot\iota$ and $\iota\ot\Delta$ to $M(A\ot A)$ and require that 
\begin{equation}
(\Delta\ot\iota)E=(E\ot 1)(1\ot E)=(1\ot E)(E\ot 1).\label{eqn:3.1a}
\end{equation}

For details, see the appendix in \cite{VD-W0}.

\opm\label{opm:3.2}
i) In the previous section, this condition would appear with $E=1$ and hence, the coproduct would be required to be non-degenerate. We didn't need this condition and in fact it followed as a consequence, see Remark \ref{opm:2.8d}.
\ssnl
ii) On the other hand, in \cite{K-VD}, we have extra assumptions for this case. The algebra $A$ is assumed to be idempotent and the  coproduct is assumed to be regular and full. As in the previous section, we do not assume these conditions and we will see how far we can get without them.  
\ssnl
iii) It is probably possible to replace these conditions with weaker ones also here, but we will not make an attempt to do so. We made an extra remark about this in Section \ref{s:concl}.
\eopm

Just as in \cite{K-VD}, we need the following extra condition on the idempotent $E$ as sitting in $M(A\ot A)$. 

\voorw\label{voorw:3.2c}
We assume the existence of two algebras $B$ and $C$, sitting in $M(A)$ in a non-degenerate way. Further we require that $E$ is a regular separability idempotent in $M(B\ot C)$ as in Definition 1.5 of \cite{VD-si.v1}. 
\evoorw

Most of the material on separability idempotents in this context is collected in  Appendix B of \cite{K-VD}. We include here the relevant properties for the convenience of the reader. More details can be found in that appendix.
\ssnl
That $B$ and $C$ sit in $M(A)$ in a non-degenerate way means that we have $BA=AB=A$ and $CA=AC=A$. It follows that the products in $B$ and $C$ are non-degenerate. The embeddings extend to embeddings of $M(B)$ and $M(C)$ in $M(A)$. So we can consider these multiplier algebras as subalgebras of $M(A)$. The same holds for $B\ot C$ and the multiplier algebra $M(B\ot C)$ is a subalgebra of $M(A\ot A)$. 
\ssnl
It is shown in \cite{K-VD} that the algebras $B$ and $C$ are determined by $E$ if they exist, see Proposition 2.1 in \cite{K-VD}.
Indeed, it can be proven that the algebras $B$ and $C$ can be defined as the left and the right leg of $E$, as sitting in $M(A\ot A)$. This is the content of the following proposition. 

\prop
If the conditions above hold, then 
\begin{align*}
&(a\ot 1)E\in A\ot C
\tussenen
E(a\ot 1)\in A\ot C\\
&E(1\ot a)\in B\ot A
\tussenen
(1\ot a)E\in B\ot A
\end{align*}
for all $a\in A$. Moreover $B$ and $C$ are the smallest subspaces of $M(A)$ for which this is true.
\eprop

From this result, we can already prove that $A$ has to be idempotent if any of the canonical maps is regular. Indeed suppose e.g.\ that $T_1$ is regular and that $\omega=0$ on $A^2$. Because $T_1$ is regular, we have 
$$\Delta(a)(p\ot q)\in A^2\ot A$$ 
and so $(\omega\ot \iota)\Delta(a)(p\ot q))=0$ for all $a,p,q$. Because we assume that $\Delta(A)(A\ot A)=E(A\ot A)$ we find that $(\omega\ot\iota)(E(p\ot q))=0$ for all $p,q$. Because $E(p\ot 1)\in A\ot C$ we already have $(\omega\ot\iota)(E(p\ot 1))=0$ for all $p\in A$. Because $E$ is full, this implies that $\omega=0$. This proves that $A=A^2$. A similar argument works for the other canonical maps.
\ssnl
It  follows from the Equations (\ref{eqn:3.1a}) that the following formulas hold.
\prop\label{prop:3.4c}
i) The subalgebras $B$ and $C$ of $A$ commute.\\
ii) For all $a\in A$ and $y\in B$ we have
$$\Delta(ay)=\Delta(a)(y\ot 1)
\tussenen
\Delta(ya)=(y\ot 1)\Delta(a).$$
iii) For all $a\in A$ and $x\in C$ we have 
$$\Delta(xa)=(1\ot x)\Delta(a)
\tussenen
\Delta(ax)=\Delta(a)(1\ot x). $$
\eprop

It is not hard to see that these properties also are true for elements in the multiplier algebras of $B$ and $C$ respectively.
\snl
We have anti-isomorphisms $S_B:B\to C$ and $S_C:C\to B$ characterized by the formulas
\begin{equation*}
E(b\ot 1)=E(1\ot S_B(b))
\tussenen
(1\ot c)E=(S_C(c)\ot 1)E.
\end{equation*}

We can define the elements $F_1, F_2, F_3$ and $F_4$ (see e.g.\ the Appendix in \cite{K-VD}). \oldcomment{Concrete reference?}{}

\defin\label{defin:3.5b}
Denote 
\begin{align}
F_1&=(\iota\ot S_C)E  \qquad\quad\text{and}\qquad\quad F_3=(\iota\ot S_B^{-1})E \label{eqn:F1andF3}\\
F_2&=(S_B\ot \iota)E \,\qquad\quad\text{and}\qquad\quad F_4=(S_C^{-1}\ot \iota)E.\label{eqn:F2andF4}
\end{align}
\edefin

These elements satisfy the following properties.

\prop\label{prop:3.5a}
The elements $F_1$ and $F_3$ belong to $M(B\ot B)$ while $F_2$ and $F_4$ belong to $M(C\ot C)$. They satisfy 
\begin{align}
E_{13}(F_1\ot 1) &= E_{13}(1\ot E) 
\tussenen 
(F_3\ot 1)E_{13}=(1\ot E)E_{13}\label{eqn:F1andF3-2}\\
(1\ot F_2)E_{13} &=(E\ot 1)E_{13}
\tussenen
E_{13}(1\ot F_4)=E_{13}(E\ot 1).\label{eqn:F2andF4-2}
\end{align}
\eprop

Finally, also observe that, with $B=C=\mathbb C 1$ and  $E=1\ot 1$, all these conditions are trivially fulfilled. This is essentially the situation in the previous section. 

\keepcomment{\nl \rood Dit vorige stuk is reeds herwerkt, maar nog niet finaal.. Wat hebben we precies nodig, wanneer en waarvoor?}{}
\nl
\bf Integrals \rm
\nl

Integrals can be defined as in Definition 2.5 of   \cite{K-VD}). However, we modify the definition slightly as we do not require the coproduct to be regular, just as we did in the previous section, see Definition \ref{defin:int}. 
\ssnl
We will use the following notations.

\notat
We denote the multiplier algebras of $B$ and $C$ by  respectively $A_s$ and $A_t$. We have seen that they are subalgebras of $M(A)$.
\enotat

Then we define integrals.

\defin\label{defin:3.7}
If $T_2$ or $T_4$ is regular, then a non-zero linear functional $\varphi$ is called a left integral if $(\iota\ot\varphi)\Delta(a)\in A_t$ for all $a\in A$. If $T_1$ or $T_3$ is regular, then a non-zero linear functional $\psi$ is called a right integral if $(\psi\ot\iota)\Delta(a)\in A_s$ for all $a\in A$.
\edefin

Remember that, given $a\in A$ and $\omega\in A'$, an element $(\iota\ot\omega)\Delta(a)$ is defined as a right multiplier if $T_2$ is regular and as a left multiplier if $T_4$ is regular. Similarly, we have that $(\omega\ot\iota)\Delta(a)$ is defined as a left multiplier if $T_1$ is regular and as a right multiplier if $T_3$ is regular. Therefore, the above definition of integrals makes sense.
\ssnl
For the integrals, defined in this way, we get the following results. They generalize the results as obtained in Proposition 2.7 in \cite{K-VD}. We use the elements $F_i$ introduced in Definition \ref{defin:3.5b}. 

\prop\label{prop:IntForm}
i) Assume that the canonical maps $T_1$ and $T_4$ are regular. If $\varphi$ is a left integral, we have 
\begin{equation*}
(\iota\ot\varphi)\Delta(a)=(\iota\ot\varphi)((1\ot a)F_4)
\end{equation*}
for all $a$ in the right leg of $\Delta$.
If $\psi$ is a right integral we get
\begin{equation*}
(\psi\ot\iota)\Delta(a)=(\psi\ot\iota)((a\ot 1)F_1)
\end{equation*}
for all $a$ in the left leg of $\Delta$.
\ssnl
ii) Assume that the canonical maps $T_2$ and $T_3$ are regular. If $\varphi$ is a left integral we have
\begin{equation*}
(\iota\ot\varphi)\Delta(a)=(\iota\ot\varphi)(F_2(1\ot a))
\end{equation*}
for all $a$ in the right leg of $\Delta$.
If $\psi$ is a right integral we get
\begin{equation*}
(\psi\ot\iota)\Delta(a)=(\psi\ot\iota)(F_3(a\ot 1))
\end{equation*}
for all $a$ in the left leg of $\Delta$.
\eprop

\bew
i) Assume that $T_1$ and $T_4$ are regular and that $\varphi$ is a left integral as in Definition \ref{defin:3.7}. Because $T_4$ is regular, given $a\in A$  we can define $y=(\iota\ot\varphi)\Delta(a)$. By definition it belongs to $A_t$ and therefore 
\begin{equation}
\Delta(yc)(p\ot q)=(y\ot 1)\Delta(c)(p\ot q)\label{eqn:3.5}
\end{equation}
 for all $c, p, q\in A$. This follows from Proposition \ref{prop:3.4c}.
 For the right hand side of Equation (\ref{eqn:3.5}) we have
 \begin{equation*}
(y\ot 1)\Delta(c)(p\ot q)
= \sum_{(a)} (a_{(1)} \ot 1)\Delta(c)(p\ot q)\,\varphi(a_{(2)}).
\end{equation*}
The element $a_{(1)}$ is covered because  the maps $T_1$ and $T_4$ are regular.
 For the left hand side of Equation (\ref{eqn:3.5}) we get
 \begin{equation*}
\sum_{(a)} \Delta(a_{(1)}c)  (p\ot q)\,\varphi(a_{(2)}).
\end{equation*}
With an extra covering factor $d$ we can write 
 \begin{equation*}
\sum_{(a)} \Delta(a_{(1)}c)  (p\ot q)\ot a_{(2)}d
=\sum_{(a)} \Delta(a_{(1)})\Delta(c)(p\ot q)\ot a_{(2)}d.
\end{equation*}
Here we use that  $T_1$ is regular. 
Now we use that $E(A\ot A)=\Delta(A)(A\ot A)$. Take $r,s\in A$ and write $E(r\ot s)=\sum_i \Delta(c_i)(p_i\ot q_i)$.
We get
\begin{align*}
\sum_{(a),i} \Delta(a_{(1)}c_i)  (p_i\ot q_i)\ot a_{(2)}d
&=\sum_{(a),i} \Delta(a_{(1)})\Delta(c_i)  (p_i\ot q_i)\ot a_{(2)}d \\
&=\sum_{(a)} \Delta(a_{(1)})E(r\ot s)\ot a_{(2)}d \\
&=\sum_{(a)} \Delta(a_{(1)})(r\ot s)\ot a_{(2)}d \\
&=\sum_{(a)} a_{(1)}r\ot  \Delta(a_{(2)})(s\ot d)
\end{align*}
We have used that $\Delta(x)E=\Delta(x)$ for all $x\in A$. For the last step we have used coassociativity with the pair $T_1,T_4$ (cf. Definition \ref{defin:1.11a}).
Because $a_{(1)}$ is covered in the two sides of this equation, we can cancel $d$ again and obtain
\begin{equation*}
\sum_{(a),i} \Delta(a_{(1)}c_i)  (p_i\ot q_i)\ot a_{(2)}
=\sum_{(a)} a_{(1)}r\ot  \Delta(a_{(2)}) (s\ot 1).
\end{equation*}
We get an element in $A\ot A\ot A$ and we can apply $\varphi$ on the last factor. We get
\begin{equation*}
\sum_{(a),i} \Delta(a_{(1)}c_i)  (p_i\ot q_i)\,\varphi(a_{(2)})
=\sum_{(a)} a_{(1)}r\ot  a_{(2)} s\,\varphi(a_{(3)}).
\end{equation*}
Putting things together, and using that  $\Delta_{13}(a)(E\ot 1)=\Delta_{13}(a)(1\ot F_4)$ (see Proposition \ref{prop:3.5a}) we find
\begin{align*}
\sum_{(a)} a_{(1)}r\ot  a_{(2)} s\,\varphi(a_{(3)})
&= \sum_{(a)} (a_{(1)} \ot 1)E(r\ot s)\,\varphi(a_{(2)})\\
&= \sum_{(a)} a_{(1)}r\ot (\iota\ot\varphi)((s\ot a_{(2)})F_4).
\end{align*}
We can apply a linear functional $\omega$ on the first factor and we get that
\begin{equation*}
 \sum_{(b)} b_{(1)} s\,\varphi(b_{(2)})
 =(\iota\ot\varphi)((s\ot b)F_4)
\end{equation*}
for all $b$ in the right leg of $\Delta$.
\ssnl
ii)  The three other cases are proven in the same way.
 \ebew
 
Remark that the formulas make sense because elements like $F_2(1\ot a)$ belong to $C\ot A$ and so $(\iota\ot\omega)(F_2(a\ot 1))$ belongs to $C$ for any linear functional $\omega$. Similarly for the other cases. In fact, we see from these formulas that not only  
$(\iota\ot\varphi)\Delta(a)$ belongs to $A_t$, the multiplier algebra of $C$, but that actually such elements already belong to $C$. This is a phenomenon we discovered already  in Section 1 in  \cite{VD-W3}. % {\blauw include a more detailed reference.} 

\opm\label{opm:3.10c}
i) In \cite{K-VD} the coproduct is assumed to be regular and full, see the beginning of Section 2 in \cite{K-VD}. Therefore we get the above results for all elements in $A$. This is used to obtain the formulas in Proposition 2.7 of \cite{K-VD}. It should be observed however that the proof given in \cite{K-VD} is not very accurate as we see from the argument in the previous proof here.
\ssnl
ii) Further we will need to have these formulas for all elements in $A$. Instead of assuming that the product is full and regular, we will just assume that the formulas hold. This boils down to giving an other definition of left and right integrals. 
\ssnl
iii) For the case $E=1$ as we have in the previous section, this just means that we define a left integral with the formula $(\iota\ot\varphi)\Delta(a)=\varphi(a)1$ and not with the requirement that $(\iota\ot\varphi)\Delta(a)$ is a scalar multiple of $1$. See Appendix \ref{sect:appA}.
\eopm

Therefore, it makes sense to \emph{work only with integrals that satisfy the formulas of the previous result for all $a\in A$}. And this is what we will do further. 
\nl
\bf The  kernels and ranges of the canonical maps \rm
\nl
We first prove the basic results about the kernels of these canonical maps. We need the concept of a left and right faithful set of linear functionals.

\defin\label{defin:3.11c}
We call a set $\mathcal F$ of linear functionals on $A$ left faithful if, given $a\in A$ we have $a=0$ if $\omega(a\,\cdot\,)=0$ for all $\omega\in\mathcal F$.  We call it right faithful if, given $a\in A$ we have $a=0$ if $\omega(\,\cdot\,a)=0$ for all $\omega\in\mathcal F$.                       
\edefin
 
\prop\label{prop:3.12d}

i) If we have a right faithful set of right integrals then $\text{Ker}(T_1)=(A\ot 1)(1-F_1)(1\ot A)$.\\
ii) If we have left faithful set of left integrals, then $\text{Ker}(T_2)=(A\ot 1)(1-F_2)(1\ot A)$.\\
iii) If we have left faithful set of right integrals then $\text{Ker}(T_3)=(1\ot A)(1-F_3)(A\ot 1)$.
iv) If we have a right faithful set of left integrals $\text{Ker}(T_4)=(1\ot A)(1-F_4)(A\ot 1)$.
\eprop

\bew
First we show that $T_1((a\ot 1)F_1(1\ot b))=T_1(a\ot b)$ for all $a,b$ in $A$. We know that $F_1(1\ot b)\in B\ot A$. We can use the Sweedler type notation $E_{(1)}\ot E_{(2)}$ for $E$. Then  $(a\ot 1)F_1(1\ot b)=aE_{(1)}\ot S_C(E_{(2)})b$ so that
\begin{align*}
T_1((a\ot 1)F_1(1\ot b))
&=\Delta(aE_{(1)})(1\ot S_C(E_{(2)})b)\\
&=\Delta(a)(1\ot E_{(1)} S_C(E_{(2)})b)=\Delta(a)(1\ot b).
\end{align*}
This proves the claim. Remark that we do not need the integrals for this result.
\ssnl
Conversely, assume that $\sum_i p_i\ot q_i\in \text{Ker}(T_1)$. This means that $\sum_i\Delta(p_i)(1\ot q_i)=0$. Multiply with $\Delta(a)$ from the left and apply a right invariant functional $\psi$ on the first leg. We find
$$\sum_i(\psi\ot\iota)(\Delta(ap_i)(1\ot q_i))=0$$
and by the formula in Proposition \ref{prop:IntForm} we find that 
$$\sum_i(\psi\ot\iota)((ap_i\ot 1)F_1(1\ot q_i))=0.$$
This holds for all right integrals $\psi$ and all $a$. Because we assume that there is a right faithful set of right integrals, it follows that
$$\sum_i(p_i\ot 1)F_1(1\ot q_i)=0.$$
This gives the right expression for the kernel of $T_1$. 
\snl
The other formulas are obtained in a similar way.
\ebew

Next we look for the ranges of the canonical maps. 

\prop\label{prop:RofT}

i) Assume that $T_1$ and $T_4$ are regular. When $\varphi$ is a left integral we have
\begin{equation*}
T_1((\iota\ot\iota\ot\varphi)(\Delta_{13}(a)\Delta_{23}(b)(1\ot c\ot 1))=E(p\ot c)
\end{equation*}
where $p=(\iota\ot\varphi)(\Delta(a)(1\ot b))$.\\
ii) Assume that $T_2$ and $T_3$ are regular. When $\psi$ is a right integral we have
\begin{equation*}
T_2((\psi\ot\iota\ot\iota)((1\ot c\ot 1)\Delta_{12}(a)\Delta_{13}(b)))=(c\ot p)E
\end{equation*}
where $p=(\psi\ot\iota)((a\ot 1)\Delta(b))$.\\
iii) Assume that $T_2$ and $T_3$ are regular. When $\varphi$ is a left integral we have
\begin{equation*}
T_3((\iota\ot\iota\ot\varphi)((1\ot c\ot 1)\Delta_{23}(a)\Delta_{13}(b)))=(q\ot c)E
\end{equation*}
where $q=(\iota\ot\varphi)((1\ot a)\Delta(b)).$\\
iv) Assume that $T_1$ and $T_4$ are regular. When $\psi$ is a right integral we have
\begin{equation*}
T_4((\psi\ot\iota\ot\iota)(\Delta_{13}(a)\Delta_{12}(b)(1\ot c\ot 1))=E(c\ot q)
\end{equation*}
where $q=(\psi\ot\iota)(\Delta(a)(b\ot 1))$.
\eprop

\bew
First assume that the maps $T_1$ and $T_4$ are regular. Then
$$\Delta_{13}(a)\Delta_{23}(b)(1\ot c\ot 1)$$
is well-defined as an element in the threefold tensor product $A\ot A\ot A$ for all $a,b,c \in A$. We can apply a left integral $\varphi$ on the last factor to get elements in $A\ot A$. On such elements we can apply the canonical map $T_1$. 
Then we find
\begin{align}
T_1((\iota\ot\iota\ot\varphi)&(\Delta_{13}(a)\Delta_{23}(b)(1\ot c\ot 1)))\\
&=(\iota\ot\iota\ot\varphi)(((\Delta\ot\iota)\Delta(a))\Delta_{23}(b)(1\ot c\ot 1))\\
&\overset{\text{(a)}}=(\iota\ot\iota\ot\varphi)((\iota\ot\Delta)(\Delta(a)(1\ot b))(1\ot c\ot 1))\\
&\overset{\text{(b)}}=(\iota\ot\iota\ot\varphi)((1\ot F_2)(\Delta_{13}(a)(1\ot c\ot b)))\\
&\overset{\text{(c)}}=(\iota\ot\iota\ot\varphi)((E\ot 1)(\Delta_{13}(a)(1\ot c\ot b)))\\
&=E((\iota\ot\varphi)(\Delta(a)(1\ot b))\ot c)\\
&=E(p\ot c)
\end{align}
with $p=(\iota\ot\varphi)(\Delta(a)(1\ot b))$. For the equality (a) we use coassociativity of the coproduct and for (c) we use Proposition \ref{prop:3.5a}. For (b) we need the formulas from Proposition \ref{prop:IntForm}. Recall that we assume these to be true for all elements $a$ in the right leg of $\Delta$, see Remark \ref{opm:3.10c}.
\ssnl
The three other formulas are proven in the same way, using the other expressions in Proposition \ref{prop:IntForm}. \oldcomment{\rood Verify!}{} 
\ebew

The start of this proof is as in the proof of Proposition \ref{prop:rTl}.
\ssnl
Recall that we use $V$ for the left leg of $\Delta$ and $W$ for the right leg of $\Delta$. We do not assume fullness of the coproduct and so, in principle, these can be proper subspaces of $A$.
\ssnl
However, with the previous result in mind, we can show the following about these subspaces. Compare this result with Proposition \ref{lem:span}. 

\prop
Suppose that there exists a set of left integrals that is left faithful. Then the span of elements of the form $(\iota\ot\varphi)(\Delta(a)(1\ot b))$ where $a,b\in A$ and where $\varphi$ is a left integral, is equal to the left leg $V$ of $\Delta$. If there is a set of left integrals that is right faithful, the span of elements $(\iota\ot\varphi)(1\ot b)\Delta(a))$ is again $V$. A similar result for right integrals will give the right leg $W$ of $\Delta$
\eprop

\bew
Assume that $\omega$ is a linear functional that kills all elements of the form $(\iota\ot\varphi)(\Delta(a)(1\ot b))$ where $a,b\in A$ and where $\varphi$ is a left integral. Replace here $b$ by $bc$. Then $\varphi(xc)=0$ for all $c$ when $x=(\omega\ot\iota)(\Delta(a)(1\ot b))$. If we assume that there is a set of left integrals that is left faithful, it follows that $x=0$. This holds for all $a,b\in A$ and consequently, $\omega$ is $0$ on $V$. This proves the first statement of the proposition. The others are proven in a completely similar way.
\ebew

Remark that we do not really need integrals. The result is true for any set of functionals with the appropriate faithfulness property.
\ssnl
If we combine this with the results of Proposition \ref{prop:RofT} we find the following. 

\prop\label{prop:3.15d}
i) Assume that $T_1$ and $T_4$ are regular. If there is a left faithful set of left integrals, 
\begin{equation*}
T_1(A\ot A)\subseteq E(V\ot A)
 \subseteq T_1(V\ot A).
\end{equation*}
ii) Assume that $T_2$ and $T_3$ are regular. If there is a right faithful set of right integrals,
\begin{equation*}
T_2(A\ot A)\subseteq  (A\ot W)E \subseteq T_2(A\ot W).
\end{equation*}
iii) Assume that $T_2$ and $T_3$ are regular. If there is a right faithful set of left integrals,
\begin{equation*}
T_3(A\ot A)\subseteq (V\ot A)E \subseteq T_3(V\ot A).
\end{equation*}
iv) Assume that $T_1$ and $T_4$ are regular. If there is a left faithful set of right integrals,
\begin{equation*}
T_4(A\ot A)\subseteq E(A\ot W) \subseteq T_4(A\ot W).
\end{equation*}
\eprop

\bew
We only consider the first case, the others are proven in a completely similar way.
\ssnl
By the very definition of the left leg of $\Delta$ we have that $T_1(a\ot b)\in V\ot A$ for all $a,b$.  And as $T_1(a\ot b)=ET_1(a\ot b)$, we also get $T_1(a\ot b)\subseteq E(V\ot A)$.
\ssnl
To show that $E(V\ot A)\subseteq T_1(V\ot A)$ we use the previous result and the first formula in Proposition \ref{prop:RofT}, together with the fact that elements of the form 
\begin{equation*}
(\iota\ot\iota\ot\varphi)(\Delta_{13}(a)\Delta_{23}(b)(1\ot c\ot 1))
\end{equation*}
belong to $V\ot A$ for all $a,b,c$. Indeed, $\Delta(b)(c\ot 1)\in A\ot A$ by the regularity of the coproduct and 
\begin{equation*}
(\iota\ot\iota\ot\varphi)(\Delta_{13}(a)(1\ot r\ot s))
\end{equation*}
belongs to $V\ot A$ for all $r,s,a$.
\ebew

If $T_1$ or $T_3$ is injective, we would get $V=A$. If $T_2$ or $T_4$ is injective, we  would get $W=A$. Then we could get our main result. However, these maps are not injective. Here we proceed as follows.

\prop\label{prop:3.16d}
i) Assume that $T_1$ and $T_4$ are regular. If there is a left faithful set of left integrals and a right faithful set of right integrals we have
\begin{equation*}
T_1(A\ot A)=E(A\ot A) \tussenen
\text{Ker}(T_1)=(A\ot 1)(1-F_1)(1\ot A).
\end{equation*}
\ssnl
ii) Assume that $T_2$ and $T_3$ are regular. If there is a left faithful set of left integrals and a right faithful set of right integrals we have
\begin{equation*}
T_2(A\ot A)=(A\ot A)E \tussenen
\text{Ker}(T_2)=(A\ot 1)(1-F_2)(1\ot A).
\end{equation*}
\ssnl
iii) Assume that $T_2$ and $T_3$ are regular. If there is a right faithful set of left integrals and a left faithful set of right integrals we have
\begin{equation*}
T_3(A\ot A)=(A\ot A)E \tussenen
\text{Ker}(T_3)=(1\ot A)(1-F_3)(A\ot 1).
\end{equation*}
\ssnl
iv) Assume that $T_1$ and $T_4$ are regular. If there is a right faithful set of left integrals and a left faithful set of right integrals we have
\begin{equation*}
T_4(A\ot A)=E(A\ot A) \tussenen
\text{Ker}(T_3)=(1\ot A)(1-F_4)(A\ot 1).
\end{equation*}
\eprop
\bew
We only prove i). The other cases can be obtained in the same way or by replacing $A$ by $A^\text{op}$ and/or $\Delta$ by $\Delta^\text{cop}$ as before.
\ssnl 
Because there is a right faithful set of right integrals, by item i) of Proposition \ref{prop:3.12d} we get the formula for the kernel of $T_1$.  
\ssnl
Because we have a left faithful set of left integrals, by item i) of Proposition \ref{prop:3.15d} we have $T_1(A\ot A)\subseteq T_1(V\ot A)$. Therefore, given $a,a'$ in $A$ we can write
\begin{equation*}
T_1(a\ot a')=\sum_i T_1(v_i\ot q_i)
\end{equation*}
with $v_i\in V$ and $q_i\in A$ for all $i$. Then we get
\begin{equation*}
(a\ot 1)F_1(1\ot a')=\sum_i(v_i\ot 1)F_1(1\ot q_i).
\end{equation*}
This is true because $p\ot q\mapsto (p\ot 1)F_1(1\ot q)$ is a projection map onto the kernel of $T_1$. Apply any linear functional $\omega$ on the second factor of this equation. Then we get $ab\in V$ where $b=(\iota\ot \omega)(F_1(1\ot a'))$. Because $E$ is assumed to be a full separability idempotent, $B$ is the span of all such elements. Therefore we have $AB\subset V$. But by assumption $AB=A$ and so $A=V$. As a consequence we get  $T_1(A\ot A)=E(A\ot A)$.
\ebew

In the first place, we now get a stronger version of the Larson-Sweedler theorem as it is found in \cite{K-VD}.

\stel
Let $A$ be a non-degenerate algebra with a regular coproduct $\Delta:A\to M(A\ot A)$. Assume that there a is regular separability idempotent $E$ as in Condition \ref{voorw:3.1c} and \ref{voorw:3.2c}. If there is a left faithful set of left integrals and a right faithful set of right integrals, the pair $(A,\Delta)$ is a weak multiplier Hopf algebra. If there is a faithful set of left integrals and a faithful set of right integrals, it is a regular weak multiplier Hopf algebra.
\estel

This result is stronger than the one in \cite{K-VD} because it is not required that the coproduct is full, neither that $A$ is idempotent.
\ssnl
There are also other possible applications. If we combine $T_1$ and $T_4$ we get the following result.

\prop
Assume that $T_1$ and $T_4$ are regular. If we have a faithful set of left integrals and a faithful set of right integrals, we have 
\begin{align*}
T_1(A\ot A)&=E(A\ot A) \tussenen
&\text{Ker}(T_1)=(A\ot 1)(1-F_1)(1\ot A)\ \\
T_4(A\ot A)&=E(A\ot A) \tussenen
&\text{Ker}(T_3)=(1\ot A)(1-F_4)(A\ot 1).
\end{align*}
\eprop

This would give a \emph{left} weak multiplier Hopf algebra. Similarly we would get a right weak multiplier Hopf algebra when $T_2$ and $T_3$ are regular. Remark however that the theory of such one-sided weak multiplier Hopf algebras has not yet been developed.
\ssnl
We have already seen in the proof of Proposition \ref{prop:3.16d} that the coproduct is full under the given conditions. In fact we also can prove the following.

\prop
Suppose that  $T_1$ is regular and  that $T_1(A\ot A)=E(A\ot A)$. Then the left leg of $\Delta$ is all of $A$.
\eprop
\bew
Assume that $(\omega\ot\iota)(\Delta(a)(1\ot b))=0$ for all $a,b\in A$. Because  $T_1(A\ot A)=E(A\ot A)$ we find that $(\omega\ot\iota)(E(p\ot q))=0$ for all $p,q\in A$.  Because $E(p\ot 1)\in A\ot C$ we already have $(\omega\ot\iota)(E(p\ot 1))=0$ for all $p\in A$. Because $E$ is full, this implies that $\omega=0$. Therefore, the left leg of $\Delta$ is all of $A$.
\ebew

Similar results hold for the other canonical maps.
\ssnl
Remark that we have already argued that $A$ has to be idempotent, see a remark preceding Proposition \ref{prop:3.4c}. 
\nl
To finish this section, we consider again the essential differences between the results of this section and the previous one. To obtain the main results in the previous section (see Theorem \ref{stel:2.8}, Proposition \ref{prop:2.9d} and Theorem \ref{stel:2.11}) we need a coproduct on a non-degenerate algebra with some regularity properties and integrals with some faithful properties. For the main result in Section \ref{s:ls-wmha} we have similar conditions, but there is an issue with the definition of the integrals. 
\ssnl
Integrals in Section \ref{s:ls-mha} are defined by a stronger requirement than in Section \ref{s:ls-wmha}. This means that either we have to assume that the coproduct is full from the very beginning, or we have to work with integrals satisfying the formulas in Proposition \ref{prop:IntForm}. 

\oldcomment{\ssnl \rood Finish with more remarks. Verwijs naar sectie 1 waar we niet moesten starten met een niet gedegeneerde $\Delta)$}{}

%%%%%%%%%%%%%%%%%%%%%
%
% Conclusions
% 
%%%%%%%%%%%%%%%%%%%%%

\section{\hspace{-17pt}. Conclusion and further research}\label{s:concl} 

Roughly speaking, a pair $(A,\Delta)$ of a unital algebra with a coproduct is a Hopf algebra if there is a faithful left and a faithful right integral. This is the Larson-Sweedler theorem as it is found in \cite{L-S}. The theorem has been generalized to multiplier Hopf algebras in \cite{VD-W-ls}, to weak Hopf algebras in \cite{B-G-L} and to weak multiplier Hopf algebras in \cite{K-VD}.
\ssnl
In this paper, we have focused on the necessary assumptions about the regularity of coproduct and the faithfulness of the integrals in order to be able to construct the counit and the antipode and thus obtaining a multiplier Hopf algebra, or more generally, a weak multiplier Hopf algebra. We get slight improvements of the known results in the two cases.
\ssnl
Among other aspects, there is the automatic fullness of the coproduct. We use a trick, used by Kustermans and Vaes in the theory of locally compact quantum groups, see \cite{Ku-Va}. 
\ssnl
In fact, as mentioned before, the treatment of locally compact quantum groups is intimately related with the Larson-Sweedler theorem in the sense that the existence of the Haar weights - the operator algebraic counterparts of the integrals - is used to develop the theory. We have explained this in greater detail in Appendix \ref{sect:appB}. 
\ssnl
Finally, let us mention some open problems related to the material of this paper. We have considered coproducts on non-degenerate algebras with various regularity properties. Up to now however there is a lack of such examples. An attempt in this direction is found in \cite{VD-refl}. Also invariant functionals in these cases, possibly with only one sided faithfulness, are not yet known. The investigation here could lead, either to peculiar examples, or to stronger results.
\ssnl
There is also the issue of non-degenerateness of the coproduct. In Section \ref{s:ls-mha} it is a consequence of the other assumptions. However, as we see from Proposition \ref{prop:A.1} in the appendix, this is so because we work with the usual definition of integrals.  In Section \ref{s:ls-wmha} we need the condition, see Condition \ref{voorw:3.1c}. One can wonder if it is possible to conclude such a condition from the other axioms as in Section \ref{s:ls-mha}. See Remark \ref{opm:3.2}.

%%%%%%%%%%%%%%%%%%%%%%%%%%%%%%%%%%%%%%%%%%%%%%%%%%%%%%%%%%%%%%%%%%%%%%
%
% Appendices 
%
%%%%%%%%%%%%%%%%%%%%%%%%%%%%%%%%%%%%%%%%%%%%%%%%%%%%%%%%%%%%%%%%%%%%%

\renewcommand{\thesection}{\Alph{section}} 

\setcounter{section}{0}

\renewenvironment{stelling}{\begin{itemize}\item[ ]\hspace{-28pt}\bf Theorem \rm }{\end{itemize}}
\renewenvironment{propositie}{\begin{itemize}\item[ ]\hspace{-28pt}\bf Proposition \rm }{\end{itemize}}
\renewenvironment{lemma}{\begin{itemize}\item[ ]\hspace{-28pt}\bf Lemma \rm }{\end{itemize}}

 \section{\hspace{-17pt}. Appendix. About integrals on algebras with a coproduct}\label{sect:appA}  

As before, we have a non-degenerate algebra $A$ and a coproduct $\Delta:A\to M(A\ot A)$. We assume that the maps $T_1$ and $T_4$ are regular. We prove the following result. 

\prop\label{prop:A.1}
Assume that $A\ot A=\Delta(A)(A\ot A)$. Let $\varphi$ be a linear functional on $A$ such that $(\iota\ot\varphi)\Delta(a)$ is a scalar multiple of the identity for all $a$. Then $(\iota\ot\varphi)\Delta(a)=\varphi(a)1$ for all $a$ in the right leg of $\Delta$.
\eprop

\bew
Let $\varphi$ be a linear functional on $A$. Assume that there is a linear functional $\lambda$ on $A$ satisfying 
\begin{equation*}
(\iota\ot\varphi)\Delta(a)(c\ot 1))=\lambda(a)c
\end{equation*}
for all $a,c$.
Apply $\Delta$  and multiply with $p\ot q$ from the right to obtain
\begin{align}
\sum_{(a)} \Delta(a_{(1)}c)(p\ot q)\,\varphi(a_{(2)})
&=\lambda(a)\Delta(c)(p\ot q)\label{eqn:2}\\
&=(\iota\ot\iota\ot\varphi)(\Delta_{13}(a)\Delta_{12}(c)(p\ot q\ot 1)).\label{eqn:3}
\end{align}
Now we use that $\Delta(A)(A\ot A)=A\ot A$. We write $r\ot s=\sum_i\Delta(c_i)(p_i\ot q_i)$. We replace $c,p,q$ by these elements and take the sum.
For the right hand side (Equation \ref{eqn:3}) we obtain
\begin{equation*}
(\iota\ot\iota\ot\varphi)(\Delta_{13}(a)(r\ot s \ot 1))=\sum_{(a)} (a_{(1)}r\ot s)\,\varphi(a_{(2)}).
\end{equation*}
We claim that for the left hand side we get 
\begin{equation*}
\sum_{(a)} (a_{(1)}r\ot a_{(2)}s)\,\varphi(a_{(3)}).
\end{equation*}
To prove the claim take another element $d$ in $A$ and write
\begin{equation*}
\sum_{(a)} \Delta(a_{(1)}c)(p\ot q)\ot a_{(2)}d
=\sum_{(a)} \Delta(a_{(1)})\Delta(c)(p\ot q)\ot a_{(2)}d.
\end{equation*}
Replace $c,p,q$ and take the sum to get
\begin{align*}
\sum_{{(a)},i} \Delta(a_{(1)}c_i)(p_i\ot q_i)\ot a_{(2)}d
&=\sum_{(a)} \Delta(a_{(1)})(r\ot s)\ot a_{(2)}d\\
&=\sum_{(a)} a_{(1)}r\ot a_{(2)}s\ot a_{(3)}d\\
&=\sum_{(a)} a_{(1)}r\ot \Delta(a_{(2)})(s\ot d).
\end{align*}
For the last step we use coassociativity. In this equation, we can cancel $d$ to get
\begin{equation*}
\sum_{{(a)},i} \Delta(a_{(1)}c_i)(p_i\ot q_i)\ot a_{(2)}
=\sum_{(a)} a_{(1)}r\ot \Delta(a_{(2)})(s\ot 1).
\end{equation*}
We can apply $\varphi$ and we see that 
\begin{equation*}
\sum_{(a)} (a_{(1)}r\ot s)\,\varphi(a_{(2)})
=\sum_{(a)} a_{(1)}r\ot (\iota\ot\varphi)(\Delta(a_{(2)})(s\ot 1)).
\end{equation*}
Finally we apply a linear functional $\omega$ on the first factor  and we replace $\sum_{(a)} \omega(a_{(1)}r) a_{(2)}$ by $b$. We get
\begin{equation*}
\varphi(b)s
=\sum_{(b)} (\iota\ot\varphi)(\Delta(b)(s\ot 1)).
\end{equation*}
This completes the proof.
\ebew

One can check that we have the necessary coverings at all places.

\opm
The conditions are natural:\\
i) A condition like $\Delta(A)(A\ot A)=A\ot A$ assures that $\Delta$ is not trivial.\\
ii) That we only get the formula for $a$ in the right leg of $\Delta$ is also natural since we only require some property of $\varphi$ on the right leg of $\Delta$. \\
iii) Indeed, if $\varphi$ is $0$ on the right leg of $\Delta$ the result can only hold if $\varphi$ is $0$.
\eopm

If the right leg of $\Delta$ is all of $A$ we get that $\varphi$ is a left integral. More precisely, we have the following consequence.

\prop
Assume that $\Delta$ is a non-degenerate homomorphism and that its legs are all of $A$. Then $\varphi$ is a left integral if and only if $(\iota\ot \varphi)\Delta(a)$ is a scalar multiplier of the identity for all $a$. Similarly $\psi$ is a right integral if and only if $(\psi\ot\iota)\Delta(a)$ is a scalar multiple of the identity for all $a$.
\eprop

In the case where a counit exists and if $A$ has an identity, the argument is simple. If we apply $\varepsilon$ on $(\iota\ot\varphi)\Delta(a)=\lambda 1$ we find immediately $\varphi(a)=\lambda$.

%%%%%%%%%%%%%%%%%%
%
% Appendix on locally compact quantum groups
%
%%%%%%%%%%%%%%%%%%%%%%

\section{\hspace{-17pt}. Appendix. Locally compact quantum groups}\label{sect:appB}  
In this appendix, we very briefly discuss the relation between the Larson-Sweedler theorem and the theory of \emph{locally compact quantum groups} as it is developed in \cite{Ku-Va}. A von Neumann algebra approach can be found in the original paper \cite{Ku-Va2} but a simpler one is found in \cite{VD-vn}. 
\ssnl
We will use the survey \cite{VD-warsaw} where references are given to the precise definitions. Here we will not go into the details.
\ssnl
We first recall the definition of a locally compact quantum group in the C$^*$-algebraic framework (see Definition 6.5 in  \cite{VD-warsaw})

\defin 
Let $A$ be a C$^*$-algebra. Consider the minimal C$^*$-tensor product $A\ot A$ and its multiplier algebra $M(A\ot A)$. A coproduct on $A$ is a non-degenerate $^*$-homomorphism from $A$ to $M(A\ot A)$ satisfying coassociativity. It is also assumed that elements of the form  $(\omega\ot\iota)\Delta(a)$ and $(\iota\ot\omega)\Delta(A)$, with $a\in A$ and $\omega$ a bounded linear functional on $A$,  are dense in $A$. Then the pair $(A,\Delta)$ is called a (C$^*$-algebraic) locally compact quantum group if there exists faithful left and right Haar weights $\varphi$ and $\psi$.
\edefin

These Haar weights are defined on the positive part of $A$ with values in $[0,\infty\,]$ satisfying additivity and certain regularity properties. The left Haar weight is left invariant and the right Haar weight is right invariant. They are assumed to be faithful. For $\varphi$ it means that $a=0$ if $a\in A$ and satisfies $\varphi(a^*a)=0$. Similarly for $\psi$.
\ssnl
If we compare these conditions with the ones we have in the Larson-Sweedler theorem for multiplier Hopf algebras, we see in the C$^*$-algebraic version of the coproduct that it is assumed to be full (in an adapted sense). In the topological framework we have the requirement that elements of the form $(\omega\ot\iota)\Delta(a)$ are dense. In the algebraic framework it is required that $A$ is spanned by elements of the form $(\omega\ot\iota)\Delta(a)$ where $\omega$ is of the form $x\mapsto f(cx)$ for $c\in A$ and $f$ any linear functional. In the algebraic framework, we need regularity of the canonical map $T_2$. In the C$^*$-algebraic setting, this is not needed because elements of the form $(\omega\ot \iota)\Delta(a)$ are automatically well-defined. Also note that  a bounded linear functional $\omega$ on $A$ is always of the form $f(c\,\cdot\,)$ for some $f$ and some $c$. 
\ssnl
Next consider the von Neumann algebraic definition of a locally compact quantum group (see Definition 6.1 in \cite{VD-warsaw}).

\defin
Let $M$ be a von Neumann algebra and consider the von Neumann tensor product $M\ot M$. A coproduct on $M$ is a normal $^*$-homomorphism $\Delta:M\to M\ot M$ satisfying coassociativity.  The pair $(M,\Delta)$ is called a (von Neumann algebraic) locally compact quantum group if there exist a faithful left and a faithful right Haar weight.
\edefin
\ssnl
The Haar  weights are assumed to be normal semi-finite weights on the von Neumann algebra and respectively left and right invariant.
\ssnl
One of the main differences between the two definitions, apart from the natural ones because we consider von Neumann algebras here, is that we do not need any such density conditions to develop the theory. In the algebraic situation, as treated in this paper, this corresponds with the fact that we do not need the fullness of the coproduct to obtain the Larson-Sweedler result. This after all is one of the main issues of this paper. 
\ssnl
Interesting to note is that the trick that makes this possible in the von Neumann algebra situation is the one used to prove that the regular representation is not just an isometric map, but actually a unitary. 
In our case, this is found in the proofs of Proposition \ref{prop:2.7c} (for multiplier Hopf algebras) and Proposition \ref{prop:3.15d} for weak multiplier Hopf algebras. 

%%%%%%%%%%%%%%%%%%%%%%%%%%%%%%%%%%%%%%%%
%
% Bestand artikel9.tex - References
%
%%%%%%%%%%%%%%%%%%%%%%%%%%%%%%%%%%%%%%%%% 

\end{document}